\documentclass[smallextended]{svjour3}       % onecolumn (second format)
\smartqed  % flush right qed marks, e.g. at end of proof
\usepackage{graphicx}
\usepackage{algorithm}
\usepackage{algorithmicx}
\usepackage{graphicx}
\usepackage{color}
\usepackage[mathscr]{euscript}
\usepackage[colorlinks,urlcolor=blue,citecolor=blue,linkcolor=blue]{hyperref}

\usepackage{makeidx}
\usepackage{mathptmx}
\usepackage{amsmath,amssymb}
\usepackage{latexsym}
\usepackage{epsfig}

\begin{document}
\title{Contingency-Constrained Unit Commitment with Post-Contingency Corrective Recourse}
%\subtitle{}
\titlerunning{Contingency-Constrained Unit Commitment with Post-Contingency Corrective Recourse}

\author{Richard Li-Yang Chen\and Neng Fan\and Ali Pinar\and Jean-Paul Watson}
\authorrunning{R.L.-Y. Chen, N. Fan, A. Pinar and J.-P. Watson}

\institute{R.L.-Y. Chen \and A. Pinar \at
Quantitative Modeling and Analysis, Sandia National Laboratories, Livermore, CA 94551\\
\email{\{rlchen,apinar\}@sandia.gov}
\and
N. Fan \at
Department of Systems and Industrial Engineering, University of Arizona, Tucson, AZ  85721\\
\email{nfan@email.arizona.edu}
\and
J.-P. Watson \at
Discrete Math and Complex Systems, Sandia National Laboratories, Albuquerque, NM 87185\\
\email{jwatson@sandia.gov}
}

\date{Received: date / Accepted: date}

\maketitle

\begin{abstract}
We consider the problem of minimizing costs in the generation unit commitment problem, a cornerstone in electric power system operations, while enforcing an $N$-$k$-$\boldsymbol \varepsilon$ reliability criterion. This reliability criterion is a generalization of the well-known $N$-$k$ criterion, and dictates that at least $(1-\varepsilon_ j)$ fraction of the total system demand must be met following the failures of $k$ or fewer system components. We refer to this problem as the Contingency-Constrained Unit Commitment problem, or CCUC. We present a mixed-integer programming formulation of the CCUC that accounts for both transmission and generation element failures. We propose novel cutting plane algorithms that avoid the need to explicitly consider an exponential number of contingencies. Computational studies are performed on several IEEE test systems and a simplified model of the Western US interconnection network, which demonstrate the effectiveness of our proposed methods relative to current state-of-the-art.
\keywords{Integer programming\and Bi-level programming\and Benders decomposition\and Unit commitment\and Contingency constraints}
\end{abstract}

\section{Introduction}
\label{sec:introduction}

Power system operations aim to optimally utilize available electricity generation resources to satisfy projected demand, at minimal cost, subject to various physical transmission and operational security constraints. Traditionally, such operations involve numerous sub-tasks, including short-term load forecasting, unit commitment, economic dispatch, voltage and frequency control, and interchange scheduling between distinct operators. Most recently, renewable generation units in the form of geographically distributed wind and solar farms have imposed the additional requirement to consider uncertain generation output, increasingly in conjunction with the deployment of advanced storage technologies such as pumped hydro. Growth in system size and the introduction of significant generation output uncertainty contribute to increased concerns regarding system vulnerability. Large-scale blackouts, such as the Northeast blackout of 2003 in North America and, more recently, the blackout of July 2012 in India, impact millions of people and result in significant economic costs. Similarly, failure to accurately account for renewables output uncertainty can lead to large-scale forced outages, as in the case of ERCOT on February 26, 2008. Such events have led to an increased focus on power systems reliability, with the goal of mitigating against failures due to both natural causes and intelligent adversaries.

Optimization methods have been applied to power system operational problems for several decades; Wood and Wollenberg \cite{Wood1996} provide a brief overview. The coupling of state-of-the-art implementations of core optimization algorithms (including simplex, barrier, and mixed-integer branch-and-cut algorithms) and current computing capabilities (e.g., inexpensive multi-core processors) enable optimal decision-making in real power systems. One notable example involves the unit commitment problem, which is used to determine the day-ahead schedule for all generators in a given operating region of an electricity grid. A solution to the unit commitment problem specifies, for each hour in the scheduling horizon (typically 24 hours), both the set of generators that are to be operating and their corresponding projected power output levels. The solution must satisfy a large number of generator (e.g., ramp rates, minimum up and down times, and capacity limits) and transmission (e.g., power flow and thermal limit) constraints, achieving a minimal total production cost while satisfying forecasted demand. The unit commitment problem has been widely studied, for over three decades. For a review of the relevant literature, we refer to \cite{Hobbs2001} and the more recent \cite{Padhy2004}. Many heuristic (e.g., genetic algorithms, tabu search, and simulated annealing) and mathematical optimization (e.g., integer programming, dynamic programming, and Lagrangian relaxation) methods have been introduced to solve the unit commitment problem. Until the early 2000s, Lagrangian relaxation methods were the dominant approach used in practice. However, mixed-integer programming implementations are either currently in use or will soon be adopted by all Independent System Operators (ISOs) in the United States to solve the unit commitment problem \cite{fercreport}.

Security constraints (i.e., which ensure system performance is sustained when certain components fail) in the context of unit commitment are now a required regulatory element of power systems operations. The North American Electric Reliability Corporation (NERC) develops and enforces standards to ensure power systems reliability in North America. Of  strongest relevance to security constraints for unit commitment is the NERC Transmission Planning Standard (TPL-001-0.1, TPL-002-0b, TPL-003-0b, TPL-004-0a, \cite{NERC2011}). The TPL specifies 4 categories of operating states, labeled A through D. Category A represents the baseline ``normal" state, during which there are no system component failures. Category B represents so called $N$-$1$ contingency states, in which a single system component has failed (out of a total of $N$ components, including generators and transmission lines). NERC requires no loss-of-load in both categories A and B, which collectively represent the vast majority of observed operational states. Categories C and D of the TPL represent more extreme states, in which multiple system components fail (near) simultaneously. Large-scale blackouts, typically caused by cascading failures, are Category D events. Such failure states are known as $N$-$k$ contingencies in the power system literature, where $k$ ($k \geq 2$) denotes the number of component failures. In contrast to categories A and B, the regulatory requirements for categories C and D are vaguely specified, e.g., ``some" loss of load is allowable, and it is permissible to exceed normal transmission line capacities by unspecified amounts for brief time periods.

The computational difficulty of security-constrained unit commitment is well-known, and is further a function of the specific TPL category that is being considered. The unit commitment problem subject to $N$-$1$ reliability constraints is, given the specific regulatory requirements imposed for category B events of the TPL, addressed by system operators worldwide. However, we observe that it is often solved approximately in practice, specifically in the context of large-scale (ISO-scale) systems \cite{personalcommunication}. For example, a subset of contingencies based on a careful engineering analysis is often used to obtain a computationally tractable unit commitment problem. Alternatively, the unit commitment problem can be solved without considering contingencies, and the solution can be subsequently ``screened" for validity under a subset of contingencies (again identified by engineering analysis). Additional constraints can then be added to the master unit commitment problem, which is then resolved; the process repeats until there is no loss-of-load in the contingency states. We raise this issue primarily to point out that even the full $N$-$1$ problem is not considered a ``solved" problem in practice, such that advances (including those introduced in this paper) in the solution of unit commitment problems subject to general $N$-$k$ reliability constraints can potentially impact the practical solution of the simpler $N$-$1$ variant.

Numerous researchers have introduced algorithms for solving both the security-constrained unit commitment problem and the simpler, related problem of security-constrained optimal power flow. In the latter, the analysis is restricted to a single time period, and binary variables relating to generation unit statuses are fixed based on a pre-computed unit commitment schedule. \cite{Capitanescu2011} provides a recent review of the literature on security-constrained optimal power flow. Of specific relevance to our work is the literature on security-constrained optimal power flow in situations where large numbers of system components fail.  This literature is mostly based on worst-case network interdiction analysis and  includes solution methods based on bi-level and mixed-integer programming (see \cite{Salmeron2004,Salmeron2009,Arroyo2010,Fan2011,Zeng2013,Zeng2014}) and graph algorithms (see \cite{Pinar2010,Bienstock2010,Fan2011,LeRoDoPi06,LePiRo08}).

Following the Northeast US blackout of 2003, significant attention was focused on developing improved solution methods for the security-constrained unit commitment problem. In particular, various researchers introduced mixed-integer programming and decomposition-based methods for more efficiently enforcing $N$-$1$ reliability, e.g., see \cite{Fu2005,Fu2006,Wang2008,Lotfjou2010,Hedman2010,ONeill2010}. However, due to its computational complexity, security-constrained unit commitment considering the full spectrum of NERC reliability standards has not attracted a comparable level of attention until very recently. Specifically, \cite{Street2011a} and \cite{Wang2012} consider the case of security-constrained unit commitment under the more general $N$-$k$ reliability criterion. Similarly, \cite{Street2011a,Street2011b} and \cite{Wang2012} use robust optimization methods for identifying worst-case $N$-$k$ contingencies.

In this paper, we extend the $N$-$k$ reliability criterion to yield the more general $N$-$k$-$\boldsymbol \varepsilon$ criterion. This new criterion dictates that for all contingencies of size $j$, $j \in \{1,\cdots, k\}$, at least $(1-\varepsilon_ j)$ fraction of the total demand must be met, with $\varepsilon_j \in [0,1]$ and $\varepsilon_1 \le \varepsilon_2 \le \cdots \le \varepsilon_k$. The primary motivation for introducing this metric is that it provides a practical  and quantifiable bound on system performance under Categories C and D TPL contingencies, and can easily be expressed in mathematical optimization models. We refer to the security-constrained unit commitment problem subject to $N$-$k$-$\boldsymbol \varepsilon$ reliability as the contingency-constrained unit commitment (CCUC) problem. In the CCUC, all contingencies with $k$ or fewer element failures (generation units or transmission lines) are implicitly considered when checking for the feasibility of post-contingency corrective recourse. The CCUC is formulated as a large-scale mixed-integer linear program (MILP). To solve the CCUC, we develop two decomposition strategies: one based on a Benders decomposition \cite{Benders1962}, and another based on cutting planes derived from the solution of power system inhibition problems  \cite{Chen2012a,Chen2014}.  We then show the computational effectiveness of our algorithms on a range of benchmark instances.

Our specific contributions, as detailed in this paper, include:
\begin{itemize}
  \item We ensure the existence of a feasible post-contingency corrective recourse, taking into consideration generator ramping constraints and the no-contingency (nominal) state economic dispatch;
  \item We consider the loss of both generation units \emph{and} transmission lines;
  \item We propose novel decomposition methods to solve the contingency-constrained CCUC efficiently, and show that models and methods proposed by \cite{Hedman2010}, \cite{ONeill2010}, \cite{Street2011a}, \cite{Street2011b}, and \cite{Wang2012} are all special cases of our general approach.
\end{itemize}

The remainder of this paper is organized as follows. In Section \ref{sec:model}, we formulate the MILP model for the contingency-constrained unit commitment problem under the $N$-$k$-$\boldsymbol \varepsilon$ reliability criterion. In Section \ref{sec:solution}, two approaches based on decomposition methods are presented for solving this large-scale MILP. In Section \ref{sec:experiments}, we test our algorithms on several IEEE test systems and a simplified model of the Western interconnection. Finally, we conclude in Section \ref{sec:conclusions} with a summary of our results and directions for future research.

%-------------------------------------------------------------------------------
%-------------------------------------------------------------------------------
\section{Problem Formulation}
\label{sec:model}
This section presents our mixed-integer linear programming model for the contingency-constrained unit commitment (CCUC) problem. In Table~\ref{tab:nom}, we introduce the core sets, parameters, and decision variables of the model. The baseline unit commitment formulation, without contingency constraints, is described in Section~\ref{sec:model-standarduc}. We discuss key concepts involving $N$-$k$-$\boldsymbol \varepsilon$ contingency analysis in Section~\ref{sec:model-contingencies}, which are subsequently illustrated on an example in Section~\ref{ieee6_example}. Finally, we combine the baseline unit commitment model with $N$-$k$-$\boldsymbol \varepsilon$ contingency analysis in Section~\ref{sec:model-ccuc}, for our contingency-constrained unit commitment model.

%\subsection{Nomenclature}
%\label{sec:model-nomenclature}
\begin{table}
\caption{Nomenclature}\label{tab:nom}
\begin{tabular}[t]{|c|l|}
\hline
\multicolumn{2}{|c|}{ \bf Sets and Indices} \\
\hline
%\begin{list}{}{\leftmargin=0em \itemindent=0em}
 $\mathcal I$&  Set of buses. Indexed by $i$ for individual buses, $i$ and $j$ for pairs of buses. \\
      $I$ &   Number of buses. $I = |\mathcal I|$. \\
$\mathcal G$ &  Set of generation units. Indexed by $g$. \\
 $G$ &  Number of generation units. $G = |\mathcal G|$. \\
$\mathcal G_i$ &  Set of generation units located at bus $i \in \mathcal I$. \\
$\mathcal E$ &  Set of directed transmission lines connecting pairs of buses. Indexed by $e$. \\
$E$ &  Number of transmission lines. $E = |\mathcal E|$. \\
$\mathcal E_{.i}$ &  Set of transmission lines oriented into bus $i \in \mathcal I$. \\
$\mathcal E_{i.}$ &  Set of transmission lines oriented out of bus $i \in \mathcal I$. \\
$i_e,j_e$ &  Tail bus $i_e$ and head bus $j_e$ of transmission line $e \in \mathcal E$. \\
$\mathcal T$ &  Set of time periods in the planning horizon. Indexed by $t \in \{1,2,\cdots,T\}$. \\
$( \mathcal I, \mathcal G, \mathcal E)$ & triple that defines a power system. \\
\hline\hline
%\end{longtable}
%\begin{longtable}{| c | l|}
\multicolumn{2}{|c|}{ \bf Parameters}\\
\hline
$B_e$, $F_e$ &    Susceptance and power flow (i.e., thermal) limit of transmission line $e$.  \\
$D_i^{t}$ &   Demand (load) at bus $i \in  \mathcal I$ at time $t $.  \\
$D^{t} =  \displaystyle \sum_{i \in \mathcal I} D_{i}^{t}$ &  Total demand, summed across all buses, in period $t $.  \\
$P^{\min}_g$, $P^{\max}_g$ &    Lower/upper limits on power output for generation unit $g $.  \\
$T_g^{d0}$, $T_g^{u0}$ &  Minimum time period generation unit $g\in  \mathcal G$ must be initially offline/online.  \\
 $T_g^d$, $T_g^u$  &   Minimum time period generation unit $g\in  \mathcal G$ must remain offline/online once  the unit \\& is shut down/started up. \\
$R_g^d,R_g^u$  & Maximum ramp-down and ramp-up rate for generation unit $g \in  \mathcal G$ between  adjacent \\& time periods.  \\
$\tilde{R}_g^d,\tilde{R}_g^u$  & Maximum shutdown/startup ramp rates for generation unit $g \in  \mathcal G$ for the time period \\&  in which $g$ is turned off/on. \\
 $C_g^u,C_g^d$ &   Fixed startup/shutdown cost for generation unit $g $.  \\
 $C_g^p(\cdot)$ &   Production cost function for generation unit $g$. \\
\hline\hline
\multicolumn{2}{|c|}{ \bf Variables}\\
\hline
  $x_g^t$ &  Binary variable indicating if a generation unit $g \in  \mathcal G$ is committed ($x_g^t=1$) or not \\& ($x_g^t=0$) at time $t $.\\
$\mathbf x^t$ &  Unit commitment decision vector for all generation units at  time $t$.\\
$\mathbf x$ &    $G\times T$ unit commitment decision vector for all generation units over all time periods. \\
 $c_g^{ut},c_g^{dt}$ & Incurred startup/shutdown cost for a generation unit $g \in  \mathcal G$ at time  $t $  (if unit $g$ is started \\& up or shut down at time $t$, the respective costs are $C_g^u$ and $C_g^d$. Otherwise, $C_g^u=C_g^d=0$.) \\
$\tilde p_g^{t}$ & No-contingency state power output by generation unit $g $  at time $t $. \\
$\tilde f_e^{t}$ & No-contingency state power  flow on transmission line $e $ at time  $t $. \\
$\tilde \theta_i^{t}$ & No-contingency state phase angle at bus $i $ at time  $t $. \\
\hline
\end{tabular}
\end{table}
%We refer to the triple $( \mathcal I, \mathcal G, \mathcal E)$ as a \emph{power system}.

\subsection{The Baseline Unit Commitment Model}
\label{sec:model-standarduc}

We now present our baseline unit commitment (BUC) formulation, without contingency constraints. Our formulation is based on the mixed-integer linear programming UC formulations introduced by \cite{Carrion2006,Wu2010,Zheng2013}. We extend these formulations to capture network transmission constraints, in the form of a DC power flow model. Our BUC model is intended to reflect steady-state operational conditions, such that the system is in a no-contingency state. Consequently, we require that the demand at each bus $i \in  \mathcal I$ must be fully satisfied, i.e., no loss-of-load is allowed.

Our BUC formulation for a power system $( \mathcal I, \mathcal G, \mathcal E)$ is given as follows:
% JPW: Reformat the below, so that the \foralls are better aligned on the right side, and separated from the constraints
% JPW: Need o introduce the f vector notation, and the theta vector notation - Q is also a functio nof x, f, p, and theta - actually, what is p below?
\begin{subequations}\label{suc}
\begin{align}
\min_{\mathbf x}\quad &\sum_{t \in \mathcal T} \sum_{g \in \mathcal G}(c_g^{ut}+c_g^{dt})+ \mathcal Q(\mathbf x) \label{suc-obj}& \\
\text{s.t.} \quad & \sum_{t=1}^{T_g^{u0}} (1-x_g^t) = 0 & \forall g \in \mathcal G  \label{suc-init-on} \\
        & \sum_{t=1}^{T_g^{d0}} x_g^t = 0  & \forall g \in  \mathcal G      \label{suc-init-off}\\
	&\sum_{t'=t}^{t+T_g^u -1} x_g^{t'}\geq T_g^u (x_g^t- x_g^{t-1}) & \forall g \in \mathcal G, t = T_g^{u0}+1, \cdots, T-T_g^u+1 \label{suc-on}\\
	&\sum_{t'=t}^T \big (x_g^{t'}-(x_g^t- x_g^{t-1}) \big)\geq 0  & \forall g \in \mathcal G, t=T-T_g^u+2, \cdots, T \label{suc-on-end}\\
	&\sum_{t'=t}^{t+T_g^d-1} (1-x_g^{t'})\geq T_g^d (x_g^{t-1} - x_g^t)  & \forall g \in \mathcal G, t=T_g^{d0}+1, \cdots, T-T_g^d+1\label{suc-off}\\
	&\sum_{t'=t}^T \big ((1-x_g^{t'})-(x_g^{t-1} - x_g^t) \big) \geq 0  & \forall g \in \mathcal G, t = T-T_g^d+2, \cdots, T \label{suc-off-end}\\                     %\\
	&c_g^{ut} \geq C_g^u (x_g^t - x_g^{t-1})  & \forall g \in  \mathcal G,t \in \mathcal  T \label{suc-c-up}\\
	&c_g^{dt} \geq C_g^d (x_g^{t-1} - x_g^{t}) & \forall g \in  \mathcal G,t \in  \mathcal T \label{suc-c-down}\\
	&c_g^{ut}, c_g^{dt} \geq 0 & \forall g \in \mathcal  G,t \in \mathcal  T  \label{buc_cost_vars}\\
	&x_g^t \in \{0,1\} & \forall g \in  \mathcal  G,t \in  \mathcal  T \label{buc_binary_const}
\end{align}
\end{subequations}

The optimization objective \eqref{suc-obj} is to minimize the sum of the startup costs $c_g^{ut}$, shutdown costs $c_g^{dt}$, and generation cost $\mathcal Q(\mathbf x)$.  Constraints \eqref{suc-init-on} - \eqref{buc_binary_const} include (in order): initial online requirements for generating units \eqref{suc-init-on}; initial offline requirements for generating units \eqref{suc-init-off}; minimum online constraints in nominal time periods for generating units \eqref{suc-on}; minimum online constraints for the last $T_g^u$ time periods \eqref{suc-on-end}; minimum offline constraints in nominal time periods for generating units \eqref{suc-off}; minimum offline constraints for the last $T_g^d$ time periods \eqref{suc-off-end}; startup cost computation \eqref{suc-c-up}; shutdown cost computation \eqref{suc-c-down}; non-negativity for startup/shutdown costs \eqref{buc_cost_vars}; and binary constraints for the on/off status of generating units \eqref{buc_binary_const}. For clarity of exposition and conciseness, we define the set $\mathcal X$ given by $\mathcal X = \{\mathbf x \in \{0,1\}^{G\times T} | \textmd{ constraints }\eqref{suc-init-on} - \eqref{buc_binary_const}\}$.

The minimum generation cost $\mathcal Q(\mathbf x)$, given a unit commitment $\mathbf x$ is constrained by a combination of DC power flow constraints and unit ramping constraints, as follows:

\begin{subequations}\label{suc-sp}
\begin{align}
\mathcal Q(\mathbf x)=\min_{{\tilde {\textbf f}}, {\tilde {\textbf p}}, {\tilde {\boldsymbol \theta}}} \quad &\sum_{g \in \mathcal G}\sum_{t \in \mathcal T} C_g^p(\tilde p_g^t)& \label{sp-obj} \\%+ \sum_{t= 1}^T\sum_{c \in \mathcal C(k)} Q(x^t,p^t, c) %\label{sp-obj}\\
\text{s.t.} \quad & \sum_{g \in   \mathcal G_i} \tilde p_g^t + \sum_{e \in   \mathcal E_{.i}}\tilde  f_{e}^t - \sum_{e \in   \mathcal E_{i.}}\tilde  f_{e}^t = D_i^t  &  \forall i \in \mathcal I, \forall t \in  \mathcal T                          \label{sp-bus-bal} \\
	& B_e (\tilde \theta_{i_e}^t - \tilde \theta_{j_e})^t - \tilde f_e^t =  0 \quad \forall e \in  \mathcal E &  \forall t \in \mathcal  T                \label{sp-flow}                                \\      %\label{sp-flow}\\
	& -F_e \leq \tilde f_e^t \leq F_e  &  \forall e \in \mathcal  E,  \forall t \in \mathcal  T    \label{sp-flow-cap}                  \\      %\label{sp-flow-cap}\\
	& P^{\min}_g x_g^t\leq \tilde p_g^t \leq P^{\max}_g x_g^t  &   \forall g \in \mathcal  G, \forall t \in \mathcal  T    \label{sp-p-bounds}   \\   %\label{sp-p-bounds}\\
	& \tilde p_g^{t} - \tilde p_g^{t-1} \leq R_g^u x_g^{t-1} + \tilde{R}_g^u (x_g^t-x_g^{t-1}) + P^{\max}_g ( 1-x_g^t) &  \forall g \in \mathcal  G, \forall t \in \mathcal  T \label{sp-ramp-up} \\    %\label{sp-ramp-up} \\
	& \tilde p_g^{t-1}- \tilde p_g^t\leq R_g^d x_g^{t}+\tilde{R}_g^d (x_g^{t-1}-x_g^t)+P^{\max}_g(1-x_g^{t-1}) &  \forall g \in  \mathcal G, \forall t  \in  \mathcal T\label{sp-ramp-down}
	%& (\mathbf f^{ct},\mathbf p^{ct},\mathbf q^{ct}, \boldsymbol \theta^{ct}) \in \mathcal Q^c(\mathbf x^t, \mathbf p^t, \mathbf d^c) \quad \forall t \in  \mathcal T, c \in TBD
\end{align}
\end{subequations}

The optimization objective \eqref{sp-obj} is to minimize generation cost given a unit commitment $\mathbf x$, where $C_g^p(\tilde p_g^t)$ is a linear approximation of  generation cost for thermal units, as is commonly employed. We discuss this linearization further below. Constraints \eqref{sp-bus-bal}-\eqref{sp-ramp-down} constitute an optimal power flow formulation, and include (in order): power balance at each bus \eqref{sp-bus-bal}; power flow on a line, proportional to the difference in voltage phase angles at the terminal buses  \eqref{sp-flow}; transmission line capacity limits \eqref{sp-flow-cap}; lower and upper bounds for committed generation unit output levels \eqref{sp-p-bounds}; and generation ramp-up/ramp-down constraints for pairs of consecutive time periods \eqref{sp-ramp-up}  and \eqref{sp-ramp-down}.

By linearizing the generation cost functions, \eqref{suc}-\eqref{suc-sp} provides a mixed-integer linear programming (MILP) formulation of the unit commitment problem with transmission constraints, but without contingency constraints. A solution to the BUC provides an on/off schedule for all generation units, over all time periods in the horizon $T$. In practice, committed generation units are adjusted on an hourly or sub-hourly basis, by ramping up or down specific units in order to satisfy realized demand. Further, additional fast-reaction (i.e., ``peaker") units can be brought online if necessary. However, this process occurs on a different time scale than the BUC, i.e., one or two hours prior to real-time execution.

\remark {Our BUC model most closely represents the reliability unit commitment problem, which ISOs and vertically integrated utilities solve  every night. In contrast, the day-ahead unit commitment problem is executed earlier in the day, and is used to clear the market and set nodal electricity prices. While there are differences between the two problem variants, specifically in terms of the inputs (e.g., bids driving aggregate demand, in contrast to ISO-forecasted load), the basic BUC model can be easily re-cast into either variant.}

\remark {The number of time periods that unit $g$ has been online/offline prior to $t=1$ should satisfy $T_g^{u0} \times T_g^{d0}=0$. That is, if a unit $g$ is online prior to time period $1$, $T_g^{u0} \ge 0$ and $T_g^{d0} = 0$. Similarly, if unit $g$ is offline before time period  $1$, $T_g^{u0} = 0$ and $T_g^{d0} \ge 0$.}

\remark {The structure of the BUC solution space is known to be degenerate, due to the nature of the phase angles $\tilde \theta_i^t$. In particular, alternative optimal solutions can be obtained by shifting all of the $\tilde \theta_i^t$ of a given optimal solution by a constant factor. To mitigate this degeneracy, and following common practice in the literature, we require in our numerical experiments that the value $\tilde \theta_r^t$ for a pre-defined ``reference'' bus $r$ be equal to 0 for $t \in \mathcal{T}$.}

\remark {Generation cost curves $C_g^p (p_g^t)$ are generally specified as quadratic functions of the form: $C_g^p (p_g^t)=c_g^{p2}(p_g^{t})^2+c_g^{p1}p_g^{t}+c_g^{p0}$. However, because the $C_g^p (p_g^t)$ are non-decreasing convex functions of $p_g^{t}$, they can be easily approximated using a piecewise linear function (see \cite{Carrion2006}). Many researchers make a further simplification by assuming a linear cost function, which corresponds to the not uncommon case in which a generator offers into the market with a single marginal cost factor. We make this assumption below in our numerical experiments, specifically that $C_g^p (p_g^{t})=c_g^p p_g^{t}$. The extension to the more general piecewise construct discussed above is straightforward, and does not impact the algorithms we introduce in Section~\ref{sec:solution}. Practically, piecewise constructs would inflate the solve times, but not significantly.}

\remark {Variables and constraints to capture reserve requirements are common in the unit commitment literature, but are absent in our unit commitment models. As noted in \cite{Hedman2010}[p. 1056], ``The primary purpose of spinning and non-spinning reserves is to ensure there is enough capacity online to survive a contingency''. Hedman et al. \cite{Hedman2010} make this argument in the context of $N$-$1$ reliabiliy; the argument for the exclusion of reserve models is even stronger for $N$-$k$ contingencies. Reserves,  specifically spinning reserves, also serve as proxies for explicitly dealing with uncertainty in demand and variable generation (e.g., wind and solar plant) output. However, again following \cite{Hedman2010}, we argue that enforcing $N$-$k$ reliability (even when $k=1$) is  likely to ensure sufficient spinning reserves are online to deal with forecast errors in both demand and variable generation. We demonstrate that this is indeed the case in Section \ref{ieee6_example} by analyzing the CCUC for a 6-bus system.}

\subsection{$N$-$k$-$\boldsymbol \varepsilon$ Contingency Constraints for Reliability Requirements}
\label{sec:model-contingencies}
According to the NERC TPL standard, in the event of a loss involving a single component (i.e., an $N$-$1$ contingency), a power system must remain stable and satisfy all demand.
% such that all thermal and voltage limits are satisfied.
%\JP{Why not generator limits? Plus, thermal limits and generation limits are allowed to increase a bit. We should perhaps also not reference voltage limits, because we don't really model that explicitly (we don't refer to voltage anywhere but at a bus, in the context of an angle.}
 %and failure to satisfy demand (i.e., loss-of-load) is not allowed.
In the case of  two or more  simultaneous losses (i.e., an $N$-$k$ contingency with $k \ge 2$), the system must maintain stability; however, a pre-planned or controlled loss-of-load is allowed. Therefore, prior to analyzing the contingency-constrained unit commitment problem, we must first ensure that the BUC model can yield solutions that satisfy such requirements.

We consider the loss of elements in a power system $(\mathcal I, \mathcal G, \mathcal E)$ in both the set $\mathcal G$ of generating units and the set $\mathcal E$ of transmission lines.
 The parameters and the variables in our formulation are defined in Table~\ref{tab:nke}.
%We introduce the following parameters to define the set of possible contingencies:
\begin{table}[thb]
\caption{ Variables and parameters $N$-$k$-$\boldsymbol \varepsilon$ contingency analysis }
\label{tab:nke}
\begin{center}
\begin{tabular}{|c|l|}
\hline
\multicolumn{2}{| l |}{\bf Parameters}\\ \hline
  $k$ &  Maximum number of simultaneous element failures. \\\hline
   $\mathcal{C}(j)$  & Set of all contingencies with \emph{exactly} $j$ failed generation units \\& and/or transmission lines  for $j \in \{1,\cdots,k\}$. Indexed by $c$.\\\hline
      $|c|$ &  Size of contingency $c$, i.e., the number of failed elements.\\\hline
  $\mathcal C = \displaystyle \cup_{j=1}^{k} \mathcal C(j)$:  & Set of all contingencies with $k$ \emph{or fewer} failed elements (generation \\& units and/or transmission lines). $|\mathcal C| = C$. \\\hline
   $\tilde d_g^c\in \{0,1\}$  & Parameter specifying whether generation unit $g \in  \mathcal G$ is involved in \\& contingency $c \in \mathcal{C}$.\\\hline
   $\tilde d_e^c\in \{0,1\}$ & Parameter specifying whether transmission line $e \in  \mathcal E$ is involved in \\& contingency $c \in \mathcal{C}$.\\\hline
   $\tilde{\textbf d}^c\in \{0,1\}^{G+T}$ &  Vector that is the concatenation of $d_g^c \ \forall g \in \mathcal G$ and $d_e^c \ \forall e \in \mathcal E$.\\\hline
    $\varepsilon_ j$  & Parameter indicating the maximum fraction of total system load that \\& can be shed in a size $j$ contingency state, for $j=1,\cdots,k$.\\\hline
     $\boldsymbol \varepsilon$  & Parameter vector indicating the maximum fraction of total  load \\& that can be shed for each contingency size, $\boldsymbol \varepsilon = (\varepsilon_ 1,\cdots,\varepsilon_ k)$. \\\hline
     $\Delta_g^j$  & Multiplicative factor applied to the ramping limits of generator $g \in \mathcal G$ \\& during a size $j \in \{1,\cdots,k\}$ contingency ($\Delta_g^j \ge 1$). \\\hline
 $\Delta_e^j$  & Multiplicative factor applied to the power flow limits of transmission \\&  line $e \in \mathcal E$ during a size $j \in \{1,\cdots,k\}$ contingency ($\Delta_e^j \ge 1$). \\
 \hline\hline
\multicolumn{2}{| l |}{\bf Variables}\\ \hline
$p_g^{ct},f_e^{ct},\theta_i^{ct}$ & corresponding values of $\tilde p_g^{t}, \tilde f_e^{t}, \tilde \theta_i^{t}$ during contingency $c \in \mathcal{C}$. \\ \hline
  $q_i^{ct}$ & Loss-of-load amount during contingency $c $ at bus $i$ at time  $t $. \\ \hline
\end{tabular}
\end{center}
\end{table}

%Variables associated with each contingency are then given as follows:

%\begin{list}{}{\leftmargin=0em \itemindent=0em}
%  \item $p_g^{ct},f_e^{ct},\theta_i^{ct}$: corresponding values of $\tilde p_g^{t}, \tilde f_e^{t}, \tilde \theta_i^{t}$ during contingency $c \in \mathcal{C}$.
%  \item $q_i^{ct}$: The loss-of-load amount during contingency $c \in \mathcal{C}$ at bus $i \in  \mathcal I$ in time period $t \in  \mathcal T$.
%\end{list}

We express the $N$-$k$ contingency set $\mathcal{C}$ as follows:
\begin{align}
\mathcal{C}=\Big \{ \tilde {\textbf{d}^c}  : \left( \sum_{g \in  \mathcal G} \tilde d_g^c + \sum_{e \in  \mathcal  E} \tilde d_e^c \right) \le k  \Big \}.
\end{align}

\remark {
% Let $N$ denote the total number of generation units and transmission lines, i.e., $N=G+E$.
The number of contingencies within the set $\mathcal{C}$ is then given by:
$${G+E\choose 1}+\cdots+{G+E\choose k}\le (G+E+1)^k-1.$$
Practically, the number of contingencies grows so rapidly  that explicit enumeration-based approaches are almost certain to fail even for modestly-sized  systems.
%For example:
%\begin{itemize}
%  \item Given $k=1$, there are $N$ contingencies in $\mathcal{C}(1)$;
% \item Given $k=2$, there are $N^2/2+N/2$ contingencies in $\mathcal{C}(2)$;
% \item Given $k=3$, there are $N^3/6+5N/6$ contingencies in $\mathcal{C}(3)$.
%\end{itemize}
}

We assume that a given contingency $c $ holds for all time periods $t \in \mathcal T$. Or alternatively, a power system must be $N$-$k$-$\boldsymbol \varepsilon$ compliant in all time periods $t \in \mathcal T$, for all contingencies $c \in \mathcal{C}$. We are not modeling specific issues relating to when a contingency may occur, how long it may last, and what corrective measures may be taken to restore functionality. Such issues can significantly expands the size and difficulty of the associated unit commitment problem, and is beyond the scope of this work. Further, generation costs are not optimized in post-contingency operation; following precedence in the literature, only constraints related to power flow on the non-contingency system elements must be enforced. In other words, the primary goal during a contingency state is operational feasibility and not cost minimization.  Additionally, multiple failure contingencies are extreme events with correspondingly low occurrence probabilities.  Therefore, consideration of the cost of these extreme events during operations planning is unnecessary, and may result in prohibitively expensive operations.

Given these assumptions, we formulate the \emph{post-contingency corrective recourse} constraints (i.e., the constraints that must be satisfied as the system state is altered in response to a contingency event, starting from a given steady state) $\mathcal R({\textbf{x}},\tilde{\textbf{p}}, \tilde {\textbf d^c})$ for a contingency prescribed by $\tilde {\textbf d^c}$, under a unit commitment decision vector $\mathbf x$ and the no-contingency state generation schedule $\tilde{\textbf{p}}$ as follows:
\begin{subequations}\label{cont_const}
\begin{align}
        \mathcal R({\textbf{x}},\tilde{\textbf{p}}, \tilde {\textbf d^c}): \quad & \sum_{g \in  \mathcal G_i} p_g^{ct} + \sum_{e \in  \mathcal E_{.i}}  f_{e}^{ct}  - \sum_{e \in  \mathcal E_{i.}}  f_{e}^{ct}  + q_i^{ct}= D_i^t \quad &  \forall i \in \mathcal I, \forall t \in  \mathcal T                          \label{csp_flow_bal} \\
        %& Gy_p^{ct} + Ay_f^{ct} +y_q^{ct}= d^t \quad \forall t \in T \label{csp_flow_bal} \\
	& B_e (\theta_{i_e}^{ct} - \theta_{j_e}^{ct} )(1- \tilde d_e^c) - f_e^{ct}=  0 \quad & \forall e \in  \mathcal E, \forall t \in  \mathcal T  \label{csp_kirchoffs} \\
	& - F_e\Delta_e^{|c|}(1-  \tilde d_e^c) \leq f_e^{ct} \leq F_e\Delta_e^{|c|}(1- \tilde d_e^c ) \quad & \forall e \in  \mathcal E, \forall t \in \mathcal  T                   \label{csp_trans_bounds} \\
	&  p_g^{ct} \leq P^{\max}_g x_g^t (1- \tilde d_g^c)\quad & \forall g \in  \mathcal G, \forall t \in  \mathcal T \label{csp_gener_bounds} \\
	&p_g^{ct}\leq R^u_g \Delta_g^{|c|} + \tilde p_g^t  \quad  & \forall g \in  \mathcal G, \forall t \in  \mathcal T \label{csp_ru}  \\
	&-p_g^{ct}\leq (R^d_g\Delta_g^{|c|} - \tilde p_g^t)(1- \tilde d_g^c) \quad & \forall g \in  \mathcal G, \forall t \in  \mathcal T    \label{csp_rd}  \\
	& q_i^{ct} \leq  D_i^t \quad & \forall i \in I, \forall t \in  \mathcal T   \label{csp_lol_i}   \\
	& \sum_{i \in  \mathcal I} q_i^{ct} \leq \varepsilon_ {|c|} D^t \quad & \forall t \in  \mathcal T \label{csp_lol_budget}\\
	& p_g^{ct} \ge 0 \quad & \forall g \in  \mathcal G, \forall t \in  \mathcal T\\
	& q_i^{ct} \ge 0 \quad & \forall i \in  \mathcal I, \forall t \in \mathcal  T.
\end{align}
\end{subequations}

Constraints \eqref{csp_flow_bal} enforce power balance at each bus, leveraging additional loss-of-load variables $q^{ct}$. Constraints \eqref{csp_kirchoffs} enforce DC power flow on those lines that are not part of contingency $\tilde {\textbf d^c}$ . Constraints \eqref{csp_trans_bounds} enforce transmission line capacity limits. If a line is not part of a contingency, then the power flow limit is given by $F_e\Delta_e^{|c|}$; otherwise, the power flow is constrained to equal zero. Constraints \eqref{csp_gener_bounds} enforce upper bounds on power output of committed (or ``on'') generation units not involved in the contingency $\tilde {\textbf d^c}$; otherwise, the power output is constrained to equal to zero. Constraints \eqref{csp_ru} enforce generator ramp-up constraints.  If a generation unit is part of the contingency, then its corresponding power output level during the contingency is zero ($p_g^{ct}=0$) and the constraint is non-binding. Otherwise, a generator can only ramp-up by $R_g^u$ from its pre-contingency level. Similarly,  constraints \eqref{csp_rd}  enforce generation ramp-down constraints.  If a generation unit is involved in a contingency, then $(1-d_g^c)=0$  and the ramp-down constraint is non-binding. Otherwise, a generator can only ramp-down by $R_g^d$ from its pre-contingency level. Constraints \eqref{csp_lol_i} ensure that loss-of-load at each bus does not exceed the demand at that bus. Finally, constraints \eqref{csp_lol_budget} ensure that at most $\varepsilon_{|c|}$ fraction of the aggregate demand can be shed in a size $|c|$ contingency.

Observe that in (\ref{cont_const}), lower limits on power output for generation units not in the contingency are relaxed to ensure sufficient operational flexibility.  These lower limits can be easily incorporated for systems with sufficiently flexible generation units.  In addition, we implicitly assume that the on/off state of generation units \emph{not} involved in a contingency are fixed and cannot be changed via recourse variables during the contingency. For those generation units that are not involved in the contingency, the power output levels $p^{ct}$ are not allowed to deviate from the baseline (pre-contingency) power output levels $\tilde p_g^{t}$ beyond the interval $[\tilde p_g^t  - R_g^d, \tilde p_g^t + R_g^u]$, given physical ramping limitations.

\subsection{A 6-Bus Illustrative Example}
\label{ieee6_example}

We now examine the impact of contingency constraints on optimal BUC solutions using the 6-bus test system introduced in \cite{Fu2005}-\cite{Fu2006}. Our goals are to concretely illustrate (a) the often significant changes in solution structure induced by the requirements to maintain $N$-$k$-$\boldsymbol \varepsilon$ in unit commitment, relative to the baseline and $N$-$1$ cases, and (b) the redundant nature of contingency constraints, in that satisfaction of one contingency state yields solutions that can ``cover" many other contingency states. The original test system consists of 6 buses, 7 transmission lines, and 3 generating units.  We modified this instance for purposes of our analysis as follows. We augmented the system with three additional, fast-ramping generators G4, G5, and G6, located at buses 1, 2, and 6, respectively. This modification ensures there is sufficient generation capacity to satisfy the $N$-$k$-$\boldsymbol \varepsilon$ criterion during contingency states. Data for the original generator set and the three additional generators is summarized in Table \ref{tab10}. Transmission line data is summarized in Table \ref{tab11}.

\begin{table}[t]

\centering
\small{
\caption{Generator data for the 6-bus test system}
\centering
\begin{tabular}{ c | c | c  c  | c c c}
\hline
Unit 		&	Bus 		 &Prod.		&Start-            	& $P^{\max}$	&  $P^{\min}$		& Ramp		 \\[-0.3ex]
		&	No. 		&(\$/MW)		&up (\$)		&(MW)		& (MW)			& (MW/h)	\\[0.2ex]
\hline
G1		&	1		&	13.51	&	125		& 220		& 100			& 55		\\[-0.3ex]
G2		&	2		&	32.63	&	249		& 100		& 10				& 50		\\[-0.3ex]
G3		&	6		&	17.69	&	0		& 100 		& 10			    	& 20		\\[-0.3ex]
G4		&	1		&	42		&	50 		& 100 		& 0			    	& 50		\\[-0.3ex]
G5		&	2		&	42		&	50 		& 100 		& 0			    	& 50		\\[-0.3ex]
G6		&	6		&	42		&	50 		& 100 		& 0			    	& 50		\\[-0.3ex]
\hline
\end{tabular}
}
\label{tab10}
\end{table}

Consistent with \cite{Fu2005}-\cite{Fu2006}, the unit shutdown cost is negligible and assumed to be zero in our analysis.  For illustrative purposes, we only consider the BUC with a single time period, with loads of 51.2, 102.4, and 42.8 at buses 3, 4, and 5, respectively. Runtime results for the full 24-hour instance are presented in Section \ref{sec:experiments}.

\begin{table}[t]
\centering
\caption{Transmission line data for the 6-bus test system}
\centering
\small{
\begin{tabular}{ c | c | c  | c  | c c c}
\hline
Line 		&	From 	&To		&$B_e$  	&$F_e$	\\[-0.3ex]
No.		&	Bus		&Bus	&		&(MW)		\\[0.2ex]
\hline
L1		&	1		&2		&5.88	& 200		\\[-0.3ex]
L2		&	1		&4		&3.88	& 100		\\[-0.3ex]
L3		&	2		&4		&5.08	& 100 		\\[-0.3ex]
L4		&	5		&6		&7.14	& 100 		\\[-0.3ex]
L5		&	3		&6		&55.56 	& 100 		\\[-0.3ex]
L6		&	2		&3		&27.03 	& 100 		\\[-0.3ex]
L7		&	4		&5		&27.03	& 100 		\\ [-0.3ex]
\hline
\end{tabular}
}
\label{tab11}
\end{table}

\begin{figure}[t]
   \centering
  % \vskip -0.5cm
    \includegraphics[width=1\textwidth, angle=0]{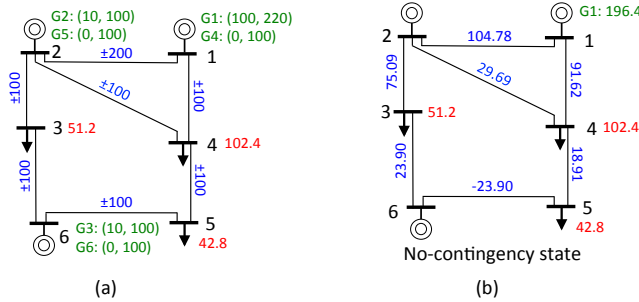}
    \vskip -4.5cm
    \caption{\footnotesize {(a) Single line diagram of the modified 6-bus test system. (b) An optimal BUC solution to the 6-bus test system, ignoring contingency constraints. Green, blue, and red text respectively refers to characteristics of generators, lines, and buses.}}
    \label{figure1}
\end{figure}

A single line diagram of the 6-bus test system is shown in Figure \ref{figure1}(a). Generator capacity bounds, transmission line capacity bounds, and loads are shown adjacent to their corresponding system elements. When contingencies are ignored, the optimal BUC solution commits a single unit (G1 at bus 1), generating 196.4 MW to meet the total demand. The no-contingency economic dispatch is shown graphically in Figure \ref{figure1}(b).

\begin{figure}[t]
   \centering
  % \vskip -0.5cm
    \includegraphics[width=1\textwidth, angle=0]{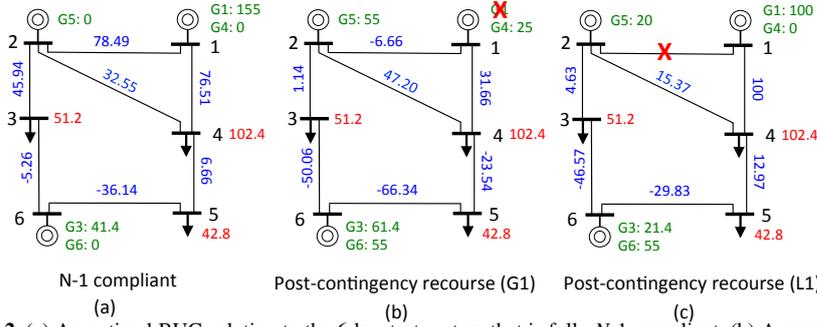}
    \vskip -4.5cm
    \caption{\footnotesize {(a) An optimal BUC solution to the 6-bus test system that is fully $N$-$1$ compliant. (b) A corrective recourse power flow after the failure of generation unit 1.  (c) A corrective recourse power flow after the failure of  the transmission line connecting buses 1 and 2.}}
    \label{figure3}
\end{figure}

In accordance with NERC's TPL standard, loss-of-load is not permitted in single-component-failure contingency states. In order for the  6-bus test system to be fully $N$-$1$ compliant, i.e., to operate the system in such a way that there exists a post-contingency corrective recourse action for \emph{all} possible $N$-1 contingencies, 5 generation units must be committed, as shown in Figure~\ref{figure3}(a). Of these, two units (G1 and G3) provide generation capacity during the no-contingency state, while three units (G4, G5, and G6) function as spinning reserves.  Unlike the practical approach of explicitly setting aside spinning reserves (e.g., equivalent to the capacity of the largest unit) via constraints, our proposed CCUC model implicitly and automatically selects units to provide spinning reserves, within the context of satisfying contingency constraints. Further, in contrast to the approach of explicitly allocating spinning reserves, our proposed CCUC model guarantees that there is adequate transmission capability to dispatch the generator outputs during all contingency states.

The optimal $N$-$1$-compliant BUC solution shown in Figure \ref{figure3}(a) represents the system in steady state operations, i.e., under no observed contingency. Figures \ref{figure3}(b) and  \ref{figure3}(c) illustrate feasible corrective recourse power flows for single-component contingency states corresponding to the failure of generation unit 1 and transmission line 1 (connecting buses 1 and 2), respectively. The total operating cost of the  $N$-$1$ compliant solution is approximately 6.52\% higher than an optimal no-contingency BUC solution.

The modified 6-bus system has 13 (7 transmission lines and 6 generators) possible single-component contingency states. We observe that it is sufficient to consider \emph{only} the two contingency states shown in Figure \ref{figure3}(b) and Figure \ref{figure3}(c) in order to achieve full $N$-$1$ compliance.  In other words, accounting for those two contingencies \emph{implicitly} yields feasible corrective recourse actions for the other $N$-$1$ contingency states. As we discuss in Section~\ref{sec:solution}, in most practical systems only a small number of contingency states are likely to impact the optimal unit commitment solution. Consequently, we design our algorithm to screen for these critical contingencies implicitly, without the need to explicitly consider all possible combinations of system component failures -- thus avoiding the combinatorial explosion in the number of possible contingencies.

\begin{figure}[t]
   \centering
  % \vskip -0.5cm
    \includegraphics[width=1\textwidth, angle=0]{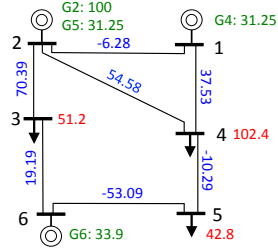}
    \vskip -5.2cm
    \caption{\footnotesize {An optimal BUC solution to the  6-bus test system that is fully $N$-$2$-$\boldsymbol \varepsilon$ compliant with $\varepsilon_ 2 = 0.27$ allowable loss-of-load.}}
    \label{figure2}
\end{figure}

If the maximum allowable contingency size is increased to $k=2$, the optimal BUC solution for the  6-bus test system commits four generation units, as shown in Figure \ref{figure2}. In addition to including $k=2$ contingencies in our analysis, we require that loads must be fully served in the no-contingency state and that a post-contingency corrective resource exist for all $k=1$ contingencies with zero loss-of-load, per TPL standards. For all $k=2$ contingencies, the allowable loss-of-load threshold is set to $\varepsilon_2 = 0.27$, to ensure that there is sufficient slack  to accommodate the loss of both transmission lines connected to bus 5. If both transmission lines connected to bus 5 fail, then the load at that bus cannot be served; the factor $0.27$ corresponds to the minimal loss-of-load under this contingency. For systems with greater redundancy and flexibility, such as those presented Section \ref{sec:experiments}, the loss-of-load threshold can be set more conservatively (i.e., lower).

Of the four committed units, one unit (G1) is producing at maximum capacity and three units (4, 5, and 6) are producing at at levels below their maximum rating. Taken together, these three units can provide up to 150MW of spinning reserves. Although fewer units are committed (4 compared to 5) relative to the $N$-$1$ solution, the two least expensive units (G1 and G2) are not committed while the three most expensive units (G4, G5, and G6) are committed in the $N$-$2$-$\boldsymbol \varepsilon$ compliant solution.

We conclude with the obvious, yet critical, observation that contingency constraints must be considered in normal (no-contingency) unit commitment operations in order to ensure that a feasible post-contingency corrective recourse exists for all contingency states under consideration.

\subsection{Contingency-Constrained Unit Commitment Formulation}\label{sec:model-ccuc}

Given the baseline unit commitment model (BUC) and associated contingency constraints as defined respectively in Sections~\ref{sec:model-standarduc} and \ref{sec:model-contingencies}, we can now describe our full contingency-constrained unit commitment (CCUC) problem:

\begin{align}\label{CCUC}
\text{CCUC}:  \quad  \min_{\mathbf x \in \mathcal X} \quad & \sum_{t \in \mathcal T} \sum_{g \in \mathcal G}(c_g^{ut}+c_g^{dt})+ \mathcal Q(\mathbf x)      \\ % JPW: The hspace is a hack - for some reason, the CCUC (7) label isn't showing up.
\text{s.t.}  \quad  & (\mathbf f, \mathbf p, \mathbf q, \boldsymbol \theta) \in \mathcal R(\mathbf x^t, {\tilde {\textbf p}}^t, \mathbf d^c) \quad \forall c \in \mathcal{C},\forall t \in \mathcal T \nonumber
\end{align}
The resulting unit commitment decision vector $\mathbf x$ represents a minimal-cost solution that satisfies (1) the non-contingency demands $D_i^t$ for each bus $i \in \mathcal I$ for each time period $t \in \mathcal T$, (2) the generation unit ramping constraints and startup/shutdown constraints, and (3) the network security and DC power flow constraints for each contingency, subject to loss-of-load allowances $\varepsilon_j$. We again note that generation costs incurred during a contingency are not considered in the objective function. Rather, only power system feasibility need be maintained, subject to the loss-of-load allowances $\varepsilon_ j$, for all $j \in \{1, \cdots, k\}$.

\section{Solution Methods}
\label{sec:solution}

The extensive formulation (EF) \eqref{CCUC} of the CCUC problem is a large-scale MILP, and has an extremely large number of variables and constraints.  For large power systems and/or non-trivial contingency budgets $k$, the EF formulation will quickly become computationally intractable. For example, the number of constraints in the second stage of the CCUC (which drives the overall problem size) is approximately given as: $C \times T(3I+2E+4G)=O\big(T\times (G+E)^k\times(I+G+E)\big)$.

Alternatively, the EF formulation of the CCUC problem has a structure that is amenable to a Benders decomposition (BD) approach, which partitions the constraints in the EF formulation into (1) a BUC problem prescribing the unit commitment decisions and the corresponding economic dispatch in the no-contingency state (this is commonly referred to as the master problem in BD), and (2) a subproblem corresponding to each contingency feasibility check, for each contingency state $c \in \mathcal C$ and time period $t \in \mathcal T$.  The BD algorithm iterates between solving the master problem (BUC), to prescribe the lowest cost unit commitment and economic dispatch, and the linear subproblems, until an optimal solution with a feasible post-contingency corrective recourse for all contingency states is obtained. In the following sub-section, we describe our Benders decomposition solution method, as it is applied to CCUC.

\subsection{A Benders Decomposition Approach}

We begin by observing that given a time period $t$, a unit commitment decision $\mathbf x^t$ and the no-contingency generation schedule $\tilde{\textbf{p}}^t$, feasibility under contingency state $c$, as prescribed by $\tilde{\textbf d^c}$, is contingent on satisfying the following DC power flow constraints.  We refer to this problem as the contingency feasibility problem {\bf{CF}}$(\mathbf x^t, \tilde{\textbf p}^t, \tilde{\mathbf d}^c)$. For conciseness of notation, we eliminate the superscript ``$ct$'' from the $f_e^{ct}, p_g^{ct}, q_i^{ct}$ and $\theta_i^{ct}$ decision variables.
\begin{subequations}\label{mod_cf}
\begin{align}
(\alpha) \quad &\sum_{g \in \mathcal G_i} p_g + \sum_{e \in \mathcal E_{.i}}  f_{e}  - \sum_{e \in \mathcal E_{i.}}  f_{e}  + q_i= D_i^t \quad &\forall i \in \mathcal  I \label{cf_bal} \\
(\beta) \quad & B_e (\theta_{i_e} - \theta_{j_e} )(1- \tilde d_e^c) -f_e=  0 \quad &\forall e \in \mathcal E \label{cf_kirk}  \\
(\hat \delta) \quad &  f_e \leq F_e \Delta_e^{|c|} (1- \tilde d_e^c ) \quad& \forall e \in \mathcal E \label{cf_f_ub}  \\
(\check \delta) \quad & -f_e \leq F_e \Delta_e^{|c|} (1- \tilde d_e^c ) \quad& \forall e \in \mathcal E \label{cf_f_lb} \\
(\gamma) \quad & p_g \leq P^{\max}_g x_g^t (1- \tilde d_g^c)\quad& \forall g \in \mathcal G  \label{cf_p_ub} \\
(\hat \lambda) \quad &p_g \leq R^u_g \Delta_g^{|c|} +  \tilde p_g^t \quad& \forall g \in \mathcal G  \label{cf_r_ub} \\
(\check \lambda) \quad &-p_g \leq (R^d_g  \Delta_g^{|c|} -  \tilde p_g^t)(1- \tilde d_g^c) \quad& \forall g \in \mathcal G  \label{cf_r_lb} \\
(\zeta) \quad &  q_i \leq  D_i^t &\quad \forall i \in \mathcal I  \label{cf_q_ub}\\
(\pi) \quad & \sum_{i \in \mathcal I} q_i \le \varepsilon_ {|c|} D^t	& \label{cf_tot_q}\\
&p_g,q_i\ge 0 &\forall g\in \mathcal G, i \in \mathcal I.
\end{align}
\end{subequations}

Using the dual variables associated with each constraint set in (\ref{cf_bal})-(\ref{cf_tot_q}), we have, by strong duality in linear programming, that $(\mathbf x^t, \tilde{\textbf p}^t, \tilde{\mathbf d}^c)$ is feasible if and only if the following dual problem  {\bf{DCF}}$(\mathbf x^t, \tilde{\textbf p}^t, \tilde{\mathbf d}^c)$ is bounded:
\begin{subequations}\label{mod_dcf}
\begin{align}
\max_{\boldsymbol \alpha,  \boldsymbol \beta,  \hat {\boldsymbol \delta}, \check {\boldsymbol \delta}, \boldsymbol  \gamma,  \hat {\boldsymbol \lambda}, \check {\boldsymbol \lambda}, \boldsymbol  \zeta, \boldsymbol  \pi} \quad &\sum_{i \in \mathcal I} D_i^t (\alpha_i+\zeta_i)+ \sum_{e \in \mathcal E} F_e\Delta_e^{|c|} (1- \tilde d_e^c)(\check \delta_e + \hat \delta_e) + \sum_{g \in \mathcal G}  P^{\max}_g x_g^t (1- \tilde d_g^c)\gamma_g   &  \\
& + \sum_{g \in \mathcal  G} \Big( (R^u_g \Delta_g^{|c|} + \tilde p_g^t)\hat \lambda_g + (R^d_g\Delta_g^{|c|}  - \tilde p_g^t )(1- \tilde d_g^c)\check \lambda_g +   \varepsilon_ {|c|} D^t \pi \Big) &\nonumber \\
\textmd{s.t.} \quad & \alpha_{i_e} - \alpha_{j_e} - \beta_e - \check \delta_e + \hat \delta_e = 0  \hspace{22ex}  \forall e \in \mathcal E  \\
& \alpha_{i_g} + \gamma_g + \hat \lambda_g - \check \lambda_g \le 0 \hspace{27ex} \forall g \in \mathcal G \\
& \alpha_i + \zeta_i \le 0  \hspace{37ex}\forall i \in \mathcal  I   \\
&\sum_{e \in \mathcal E_{i.}} B_e(1- \tilde d_e^c)\beta_e - \sum_{e \in \mathcal E_{.i}}B_e(1- \tilde d_e^c)\beta_e=0  \hspace{10ex} \forall i \in  \mathcal   I\\
%&\delta_e^\ell, \delta_e^u, \gamma_g^\ell, \gamma_g^u, \lambda_g^u, \lambda_g^d, \zeta_i, \eta_g \le 0 & \forall e \in E, g \in G, i \in \mathcal N.
&\hat {\boldsymbol \delta}, \check {\boldsymbol \delta}, \boldsymbol  \gamma,  \hat {\boldsymbol \lambda}, \check {\boldsymbol \lambda}, \boldsymbol  \zeta, \boldsymbol  \pi\le 0
\end{align}
\end{subequations}

Note that the feasible domain for  {\bf{DCF}}$(\mathbf x^t, \tilde{\textbf p}^t, \tilde{\mathbf d}^c)$, is a polyhedral cone and any solution in the domain is a ray.  By Minkowski's theorem, every such ray can be expressed as a non-negative linear combination of the extreme rays of the domain. Therefore, the dual problem {\bf{DCF}}$(\mathbf x^t, \tilde{\textbf p}^t, \tilde{\mathbf d}^c)$ is bounded if and only if  its optimal objective value is less than or equal to zero. And this happens   if and only if
\begin{align}\label{f_cut}
& \sum_{i \in \mathcal I} D_i^t (\alpha_i+\zeta_i)+ \sum_{e \in \mathcal E} F_e\Delta_e^{|c|} (1- \tilde d_e^c)(\check \delta_e + \hat \delta_e) + \sum_{g \in \mathcal G}  P^{\max}_g x_g^t (1- \tilde d_g^c)\gamma_g     \\
& + \sum_{g \in \mathcal G} \Big( (R^u_g \Delta_g^{|c|} + \tilde p_g^t)\hat \lambda_g + (R^d_g\Delta_g^{|c|}  - \tilde p_g^t )(1- \tilde d_g^c)\check \lambda_g +   \varepsilon_ {|c|} D^t \pi \Big) \le 0. \nonumber
\end{align}

We call these the Benders feasibility cuts or $f$-$cut$ for short. Below we outline the Benders decomposition algorithm as it applied to CCUC.
\begin{algorithm}%[H]
\caption{\emph{Benders Decomposition Algorithm (BD)}}
\begin{algorithmic}[1]
\State  Initialization: let $\ell \leftarrow 1$
\State Solve BUC$_\ell$
\State \textbf{if} BUC$_\ell$ is infeasible, CCUC has no feasible solution, EXIT
\State  \textbf{else}, let  $(\mathbf{x}_\ell, \tilde{\textbf{p}}_\ell)$ be an optimal solution of BUC$_\ell$
\State \hskip 0.8cm \textbf{for} each $c\in \mathcal{C}$, $t \in \mathcal T$,
\State \hskip 1.6cm solve  {CF${(\mathbf{x}^t_\ell, \tilde{\textbf{p}}^t_\ell,\tilde{\mathbf{d}}^c)}$}   and let $w^*$ be the optimal objective value
\State \hskip 1.6cm \textbf{if}  $w^*$ unbounded  add $f$-cut \eqref{f_cut} to \text{BUC$_\ell$}
\State \hskip 1.6cm \textbf{end if}
\State \hskip 0.8cm \textbf{end for}
\State \hskip 0.8cm \textbf{if} $f$-cut(s) added in (7), let $\ell \leftarrow \ell+1$ and return to (2)
\State \hskip 0.8cm \textbf{else}, $(\mathbf{x}_\ell, \tilde{\textbf{p}}^t_\ell)$ is an optimal solution, EXIT
\State \hskip 0.8cm  \textbf{end if}
\State \textbf{end if}
\end{algorithmic}
\end{algorithm}

\subsection{A Cutting Plane Method Based on the Power System Inhibition Problem}
Even with a BD approach CCUC  may not be tractable for practical size power systems, since  for every contingency $c\in \mathcal{C}$ and time period $t \in \mathcal T$, we need to ensure that a feasible DC power flow with limited loss-of-load exists, which is intractable in most cases.
% it is intractable to check power flow feasibility  for every possible contingency state $c \in \mathcal C$ and time period $t \in \mathcal T$ pair.

 In this section, we describe a cutting plane algorithm that uses a bi-level separation problem to {\em implicitly} identify a contingency state that would result in the worst-case loss-of-load for each contingency size $j$ , $j \in \{1,\cdots,k\}$.  If the worst-case generation shedding is non-zero and/or loss-of-load is above the given contingency budget $\varepsilon_j$, the current solution is infeasible, and we generate a cutting plane, corresponding to $f$-$cut$ \eqref{f_cut} to add to the BUC to protect against this particular contingency state.

\subsubsection{The Bi-Level Power System Inhibition Problem}
Given a time period $t \in \mathcal T$ and a contingency budget $j \in \{1,\cdots,k\}$, unit commitment $\mathbf x^t$, and the no-contingency generation levels $\tilde{\textbf{p}}^t$, we solve a bi-level \emph{power system inhibition problem} (PSIP), to determine the worst-case generation/load  shedding under a contingency with \emph{exactly} $j$ failed elements. In the context of PSIP, the contingency vector $\mathbf d$ is no longer a parameter but a vector of upper-level decision variables. In PSIP, the upper-level decisions $(\mathbf d)$ correspond to binary contingency selection decisions and the lower level decisions $(\mathbf f, \mathbf p, \mathbf q, \mathbf r, s, \boldsymbol \theta)$ correspond to recourse generation schedule and DC power flow under the state prescribed by the the unit commitment decisions $\mathbf x^t$, the no-contingency state economic dispatch $(\tilde{\textbf p}^t)$, and upper-level contingency selection variables $(\mathbf d)$.

Before we introduce the power system inhibition problem (PSIP), we augment the direct current power flow constraints as follows.  We introduce three sets of non-negative, continuous variables corresponding to generation shedding $r_g$ for all $g \in \mathcal G$, loss-of-load at each bus $q_i$ for all $i \in \mathcal I$ and auxiliary variable $s$ corresponding to total system loss-of-load above the allowable threshold $\varepsilon_ {j}D^t$. These variables in conjunction with additional constraints ensure that PSIP has a feasible recourse power flow for any unit commitment $\mathbf x^t$, no-contingency state economic dispatch ${\tilde{\textbf p}^t}$ and upper-level contingency selection decisions $\mathbf d$. We now state the power system inhibition problem.
\begin{subequations} \label{bpsip}
\begin{align}
\textrm{B-PSIP}({\textbf{x}}^t, \tilde{\textbf{p}}^t, j): &\nonumber  \\
\max_{\mathbf d} \ \min_{\mathbf f,\mathbf p, \mathbf q, \mathbf r, s, \boldsymbol \theta} \quad &  \sum_{g \in \mathcal G} r_g + s \label{bpsip_obj}  \\
\text{s.t.}  \quad
        & \sum_{e \in\mathcal  E} d_e + \sum_{g \in \mathcal G} d_g = j\label{bpsip_budget} \\
    	&\sum_{g \in \mathcal G_i} (p_g-r_g) + \sum_{e \in \mathcal E_{.i}}  f_{e} - \sum_{e \in \mathcal E_{i.}}  f_{e} + q_i= D_i^t  &\forall i \in \mathcal  I  \label{bpsip_bal} \\
    	& B_e (\theta_{i_e} - \theta_{j_e})(1-d_e)-f_e=  0  & \forall e \in \mathcal E  \label{bpsip_vol} \\
	& -f_e \leq F_e \Delta_e^j (1-d_e)  & \forall e \in \mathcal E  \label{bpsip_f_lb} \\
	&  f_e \leq F_e(1-d_e) \Delta_e^j   & \forall e \in \mathcal E  \label{bpsip_f_ub} \\
	& p_g \leq P^{\max}_g   x_g^t (1-d_g)  & \forall g \in \mathcal G \label{bpsip_p_ub} \\
	&p_g \leq R^u_g \Delta_g^j  + \tilde p_g^t   & \forall g \in \mathcal G \label{bpsip_ru} \\
	&-p_g \leq R^d_g \Delta_g^j  - \tilde p_g^t(1-d_g)   & \forall g \in \mathcal G \label{bpsip_rd} \\
	&  q_i\leq  D_i^t   & \forall i \in \mathcal I \label{bpsip_q} \\
	& r_g - p_g \leq  0  & \forall i \in \mathcal I \label{bpsip_r} \\
 	& \sum_{i \in \mathcal I} q_i - s \le \varepsilon_ {j} D^t \label{bpsip_s}\\
	& p_g \ge 0,\; q_i \ge 0,\; r_g \ge 0,\; s \ge 0  & \forall i \in \mathcal I, g \in \mathcal G \label{bpsip_non_neg}\\
	& d_e \in \{0,1\}, \;d_g \in \{0,1\} \quad & \forall e \in \mathcal E, \forall g \in \mathcal G. \label{bpsip_d}
\end{align}
\end{subequations}

The bi-level objective \eqref{bpsip_obj} seeks to maximize, the minimum generation shedding and loss-of-load quantity above the allowable threshold. Since $r_g$ for all $g\in \mathcal G$ and $s$ are non-negative variables, the objective value is bounded below by zero.  If the objective value is equal to zero, the current solution $({\textbf{x}}^t, \tilde{\textbf{p}}^t)$ has a feasible corrective recourse for all contingencies of size $j$.  Otherwise, the current solution cannot survive the contingency prescribed by upper-level contingency selection variables $\mathbf d$. Given a contingency state defined by $\mathbf d$, the objective of the power system operator (the inner minimization problem) is to find a corrective recourse power flow such that generation shedding and loss-of-load quantity above the allowable threshold is minimized. \eqref{bpsip_budget} is a budget constraint on the number of power system elements, generation and/or transmission, that must be in the selected contingency state. Constraints \eqref{bpsip_bal} enforce power balance at each bus, with additional generation shedding variables $r_g \in \mathcal G$ for each generator located at a bus and a bus load-shedding variable $q_i \in \mathcal I$ to ensure system feasibility.  Constraints \eqref{bpsip_vol}-\eqref{bpsip_q} are as stated in (\ref{cont_const}).  Constraints (\ref{bpsip_r}) restrict the  amount of generation shedding to be at most the generation output for each generator $g \in \mathcal G$.  Constraint (\ref{bpsip_s}) defines the amount of load shedding above the allowable threshold. If  $\sum_{i \in \mathcal I} q_i > \varepsilon_ {|c|} D^t$ then $ s = \sum_{i \in \mathcal I} q_i - \varepsilon_ {|c|} D^t$, otherwise, $s=0$.

%Remember that a failed element was assumed to be fail throughout all $t \in \mathcal T$ in our CCUC model. Therefore,  a set of generating units and/or transmission lines is selected to be failed for a contingency should consider all $t$.

\remark {Observe that  \eqref{bpsip_vol} are nonlinear constraints with terms associated with products of binary contingency-selection variables (upper level) and continuous voltage phase angles  variables (lower level); thus, B-PSIP is a bi-level, nonlinear program. }

\remark {Observe that B-PSIP is feasible for any first-stage solution $(\mathbf x^t, \tilde{\textbf{p}}^t)$ and any contingency prescribed by $\mathbf d$; the  solution $\mathbf f = \mathbf 0, \mathbf p^t = \mathbf r^t = \tilde{\textbf{p}}^t , \mathbf q^t = D^t$,  $s=(1-\varepsilon_j)D^t$ and $\boldsymbol \theta = \mathbf 0$ is feasible under any contingency state. }

Bi-level programs, such as (\ref{bpsip}), cannot be solved directly. Next, we describe a reformulation strategy to derive a mixed-integer linear programming equivalent for B-PSIP. We begin by fixing the upper level variables $\mathbf d$ and dualizing the inner minimization problem.  This results in a single-level, bilinear program with bilinear terms in the objective function.  In the resulting reformulation, there are five nonlinear terms, which are products of
binary contingency selection variables $(d_e, d_g)$ and continuous dual variables $(\beta, \hat \delta, \check \delta, \gamma, \check \lambda)$. Each of these non-linear terms can be linearized using the following strategy.

Let $u$ and $v$ be two continuous variables and $b \in \{0,1\}$.  Then the bilinear term, $(1-b)u$, can be linearized as follows. Letting $v = (1-b)u$, we introduce the following four constraints to linearize the bilinear term $(1-b)u$.
\begin{subequations}
\begin{align}
	  u-Ub \;\leq\;& v \;\le\; u + Ub \label{lin13}\\
	 -U(1-b) \;\leq\; &  v \;\le\; U(1-b) \label{lin24}
%	& v \le  \label{lin3}\\
%	& v \le U(1-b) \label{lin4}
\end{align}
\end{subequations}

\noindent
Here, parameter $U$ represents an upper bound for continuous variable $u$ and satisfies $U \ge |u|$. Assessing these  constraints for both binary values of $b$ shows that they provide a linearization.  If $b=1$, then constraint  (\ref{lin24})  implies that $v=0$.  With $v=0$, constraints  (\ref{lin13}) implies that $-U \le u\le U$, which are never binding. If $b=0$, then constraints (\ref{lin13}) implies $u=v$ and constraint (\ref{lin24})   implies $-U \le v \le U$, which are never binding.

\remark {If the bilinear term is a product of a binary variable $b$ and a non-positive variable $u$, i.e. $u \le 0$, the lower bound in  (\ref{lin24}) is redundant, and thus, can be eliminated. Analogously, if $u$ is a non-negative variable, i.e. $u \ge 0$, the upper bound in  (\ref{lin24}) is redundant, and thus, can be eliminated.}

We follow a similar strategy to linearize all five bilinear terms $(\beta, \hat \delta, \check \delta, \gamma, \check \lambda)$.  Define continuous variables $(r^1, r^1, r^3, r^4, r^5)$ and let $r_e^1 = (1-d_e) \beta_e$, $r_e^2 = (1-d_e) \hat \delta_e$, $r_e^3 = (1-d_e) \check \delta_e$, $r_g^4 =(1-d_g) \gamma_g$ and $ r_g^5 = (1-d_g) \check \lambda_g$.  Following the same linearization strategy introduced above, we now state the full mixed-integer linear PSIP formulation for completeness.
\begin{subequations}\label{mod_psip_full}
\begin{align}
\!\!\!\!\! \!\!\!\!\! \!\!\!\!\! \!\!\!\!\!
\max_{\mathbf d, \boldsymbol \alpha,  \boldsymbol \beta,  \hat {\boldsymbol \delta}, \check {\boldsymbol \delta}, \boldsymbol  \gamma,  \hat {\boldsymbol \lambda}, \check {\boldsymbol \lambda}, \boldsymbol  \zeta, \boldsymbol  \pi} \quad &\sum_{i \in \mathcal I} D_i^t (\alpha_i+\zeta_i)+ \sum_{e \in \mathcal E} F_e\Delta_e^{j} (r_e^2 + r_e^3)  + \sum_{g \in  \mathcal G}  P^{\max}_g x_g^t r_g^4\gamma_g \!\!&   \\
& + \sum_{g \in  \mathcal G} \Big( (R^u_g \Delta_g^{j} + \tilde p_g^t)\hat \lambda_g + (R^d_g\Delta_g^{j}  - \tilde p_g^t ) r_g^5 +   \varepsilon_ {j} D^t \pi \Big) &\nonumber \\
\textmd{s.t.} \quad & \sum_{e\in  \mathcal E} d_e + \sum_{g\in  \mathcal G} d_g = j & \label{psip_full_budget} \\
& \alpha_{i_e} - \alpha_{j_e} - \beta_e - \check \delta_e + \hat \delta_e = 0   & \forall e \in \mathcal E  \\
& \alpha_{i_g} + \gamma_g + \hat \lambda_g - \check \lambda_g \le 0 & \forall g \in \mathcal G \\
& \alpha_i + \zeta_i \le 0 & \forall i \in \mathcal  I   \\
& -\alpha_{i_g} + \eta_g \le 1 & \forall g \in \mathcal G\\
& -\pi \le 1 & \\
&\sum_{e \in  \mathcal E_{i.}}B_e r_e^1 - \sum_{e \in  \mathcal E_{.i}}B_e r_e^1 = 0 & \forall i \in  \mathcal   I \\
& r_e^1  \ge \max\{ \beta_e - U d_e, - U (1-d_e)   \} & \forall e \in \mathcal E \\
%& r_e^1  \ge \beta_e - U d_e & \forall e \in \mathcal E \\
%&r_e^1  \ge - U (1-d_e) & \forall e \in \mathcal E \\
 & r_e^1  \le \min \{ \beta_e + U d_e,  U (1-d_e)\} & \forall e \in \mathcal E \\
%&r_e^1  \le U (1-d_e) & \forall e \in \mathcal E \\
&r_e^2  \ge \max \{ \hat \delta_e - U d_e, - U (1-d_e)  \} & \forall e \in \mathcal E \\
%&r_e^2  \ge - U (1-d_e) & \forall e \in \mathcal E \\
&r_e^2  \le \hat \delta_e + U d_e& \forall e \in \mathcal E \\
%&r_e^2  \le U (1-d_e) & \forall e \in \mathcal E \\
&r_e^3  \ge \max\{ \check \delta_e - U d_e, - U (1-d_e) \}& \forall e \in \mathcal E \\
%&r_e^3  \ge - U (1-d_e) & \forall e \in \mathcal E \\
&r_e^3  \le \check \delta_e + U d_e& \forall e \in \mathcal E \\
%&r_g^3  \le U (1-d_e) & \forall e \in \mathcal E \\
&r_g^4  \ge \max \{ \gamma_g - U d_g, - U (1-d_g) \} & \forall g \in \mathcal G \\
%&r_g^4  \ge - U (1-d_g) & \forall g \in \mathcal G \\
&r_g^4  \le\gamma_g + U d_g& \forall g \in \mathcal G \\
%&r_g^4  \le U (1-d_g) & \forall g \in \mathcal G \\
&r_g^5  \ge \max \{ \check \lambda - U d_g,  - U (1-d_g)\} & \forall g \in \mathcal G \\
%&r_g^5  \ge - U (1-d_g) & \forall g \in \mathcal G \\
&r_g^5  \le \check \lambda + U d_g& \forall g \in \mathcal G\\
%&r_g^5  \le U (1-d_g) & \forall g \in \mathcal G
& \hat {\boldsymbol \delta}, \check {\boldsymbol \delta}, \boldsymbol  \gamma,  \hat {\boldsymbol \lambda}, \check {\boldsymbol \lambda}, \boldsymbol  \zeta, \boldsymbol  \pi\le 0
\end{align}
\end{subequations}

Next, we outline an algorithm for \emph{optimally} solving problem CCUC that combines a Benders decomposition with the aid of an oracle given by \eqref{mod_psip_full}, which acts as a separation problem.  A given solution $({\textbf{x}}^t, \tilde{\textbf{p}}^t)$ is feasible if the oracle cannot find a contingency of size $j$ that results in a loss-of-load above the allowable threshold. That is, if the optimal objective value is zero.  For each contingency budget $(j=1,\cdots,k)$, we can check for the worst-case $j$-element contingencies by solving $(\ref{mod_psip_full})$ using a failure budget of $j$ (i.e. the right-hand side of inequality (\ref{psip_full_budget}) is set to $j$, as it is right now).  Whenever the oracle determines that the current solution is \emph{not} $N$-$k$-$\boldsymbol \varepsilon$ compliant,  it returns a contingency state, prescribed by $\mathbf d$, that results in a generation shedding and/or  loss-of-load, above the allowable threshold $\varepsilon_ j D^t$ for $j$-element failures.

\subsubsection{Contingency Screening Algorithms}
We now present a cutting plane algorithm,  referred to as the \emph{Contingency Screening Algorithm 1 (CSA1)} to solve CCUC  implicitly by screening for the worst-case contingency, in terms of total generation and load shedding.
\begin{algorithm}%[H]
\caption{\emph{Contingency Screening Algorithm 1 (CSA1)}}
\begin{algorithmic}[1]
\State  Initialization: let $\ell \leftarrow 1$
\State Solve BUC$_\ell$
\State \textbf{if} BUC$_\ell$ is infeasible, CCUC has no feasible solution, EXIT
\State  \textbf{else}, let  $( {\textbf{x}}_\ell, \tilde{\textbf{p}}_\ell)$ be an optimal solution of BUC$_\ell$
\State \hskip 0.8cm \textbf{for} all $j \in \{1,\cdots, k\}$, $t \in \mathcal T$,
\State \hskip 1.6cm solve  {PSIP${( {\textbf{x}}^t_\ell, \tilde{\textbf{p}}^t_\ell,j)}$}   and let $w^*$ be the optimal objective value
\State \hskip 1.6cm \textbf{if}  $w^* > 0,$
\State \hskip 2.4cm Add $f$-cut \eqref{f_cut} to \text{BUC$_\ell$}
\State \hskip 1.6cm \textbf{end if}
\State \hskip 0.8cm \textbf{end for}
\State \hskip 0.8cm \textbf{if} $f$-cut(s) added in (7), let $\ell \leftarrow \ell+1$ and return to (2)
\State \hskip 0.8cm \textbf{else}, $({\textbf{x}}_\ell, \tilde{\textbf{p}}_\ell)$ is an optimal $N$-$k$-$\boldsymbol \varepsilon$ compliant solution, EXIT
\State \hskip 0.8cm  \textbf{end if}
\State \textbf{end if}
\end{algorithmic}
\end{algorithm}

\subsubsection{Contingency Sharing Using a Dynamic Contingency List}

\begin{figure}[t]
   \centering
  \vskip -0.5cm
    \includegraphics[width=0.7\textwidth, angle=0]{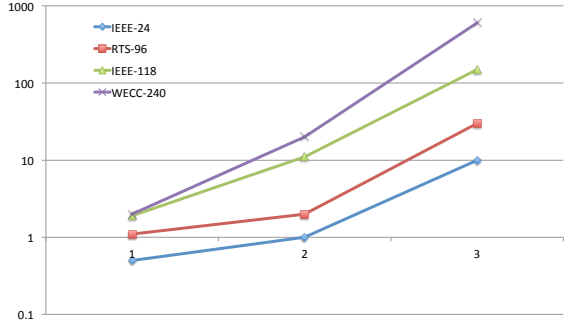}
    \vskip -1cm
    \caption{\footnotesize {Average run times (sec.) of PSIP for varies power systems and contingency budgets ($k=1,2,3$).}}
    \label{figure0}
\end{figure}

In preliminary testing using CSA1, we found that run time is significantly impacted by the need to solve a large number of PSIP instances at each master iteration of the algorithm. Specifically, we solve one instance of PSIP for each contingency-size and period pair $(j,t)$ for each master iteration.  The solution time of PSIP, as expected, is heavily impacted by the size of the power system $(\mathcal I, \mathcal G, \mathcal E)$.  Figure \ref{figure0} shows the solution time (on a logarithmic scale) of PSIP for various power system sizes and maximum contingency budgets $k$.

In solving CCUC using CSA1 we also made three observations.  First, the majority of the the total run-time was spent solving PSIP (\ref{mod_psip_full}) instances.  Secondly, a contingency $c$ that fails the system in one time period $t$ often fails the system in other time times as well, which suggests sharing of contingencies across time periods.  Thirdly, in the final solution only a small number of contingencies are actually identified.  That is to say, it is often prudent to consider a small number of contingencies \emph{explicitly} in solving CCUC.
Based on these observations, we found that it is most efficient to develop a version of the CSA algorithm that minimizes the number times we solve PSIP (\ref{mod_psip_full}) instances and allows for sharing of contingencies across time periods.  We achieve this buy using a dynamic contingency list.

We begin with an empty contingency list $\mathscr L$.  At each master iteration, we first screen all contingencies in the list for each time period $t \in \mathcal T$.  For each time period $t$, we generate feasibility cuts  \eqref{f_cut} for each violated contingency in the list.  If none of the contingencies in the list is violated in any time period $t$, we proceed to solving PSIP instances to identify a new  violated contingency. This simple procedure ensures that each violated contingency identified by solving PSIP is never redundant.  That is to say, the new contingency is not in our existing contingency list.  When a new contingency is identified, we add it to the contingency list and check for its violation in all other time periods by solving a linear DCF problem.
Our computational results indicate that this procedure results in the fewest total PSIP instances solved on average, which results in the fastest run
time. The key idea is that this procedure avoids redundant PSIP solutions to re-identify violated contingencies. This algorithm is referred to as the {\em Contingency Screening Algorithm 2 (CSA2)}.

\begin{algorithm}%[H]
\caption{\emph{Contingency Screening Algorithm 2 (CSA2)}}
\begin{algorithmic}[1]
\State  Initialization: $\ell \leftarrow 1$, $\mathscr L = \emptyset$
\State Solve BUC$_\ell$
\State \textbf{if} BUC$_\ell$ is infeasible, CCUC has no feasible solution, EXIT
\State  \textbf{else}, let  $( {\textbf{x}}_\ell, \tilde{\textbf{p}}_\ell)$ be an optimal solution of BUC$_\ell$
\State \hskip 0.8cm \textbf{for} each $c\in \mathcal{C}$, $t \in \mathcal T$,
\State \hskip 1.6cm solve  {CF${(\mathbf{x}^t_\ell, \tilde{\textbf{p}}^t_\ell,\mathbf{d}^c)}$}   and let $w^*$ be the optimal objective value
\State \hskip 1.6cm \textbf{if}  $w^* > 0,$
\State \hskip 2.4cm Add $f$-cut \eqref{f_cut} to \text{BUC$_\ell$}
\State \hskip 1.6cm \textbf{end if}
\State \hskip 0.8cm \textbf{end for}
\State \hskip 0.8cm \textbf{if} $f$-cut(s) added in step (7)
\State \hskip 1.6cm let $\ell \leftarrow \ell+1$, return to step (2)
\State \hskip 0.8cm \textbf{end if}
\State \hskip 0.8cm \textbf{for} all $j \in \{1,\cdots, k\}$, $t \in \mathcal T$,
\State \hskip 1.6cm solve  {PSIP${( {\textbf{x}}^t_\ell, \tilde{\textbf{p}}^t_\ell,j)}$}   and let $z^*$ be the optimal objective value
\State \hskip 1.6cm \textbf{if}  $z^* > 0,$  add $f$-cut \eqref{f_cut} to \text{BUC$_\ell$}
\State \hskip 2.4cm let $\ell \leftarrow \ell+1$,  $\mathscr L \leftarrow \mathscr L \cup \{c\}$, return to step (2)
\State \hskip 1.6cm \textbf{end if}
\State \hskip 0.8cm \textbf{end for}
\State  \textbf{end if}
\State $( \tilde{\textbf{x}}_\ell, \tilde{\textbf{p}}_\ell)$ is an optimal $N$-$k$-$\boldsymbol \varepsilon$ compliant solution, EXIT
\end{algorithmic}
\end{algorithm}

\section{Computational Experiments}
\label{sec:experiments}

We implemented our proposed models and algorithms in C++ using IBM's Concert Technology Library 2.9 and CPLEX 12.1 MILP solver. All experiments were performed on a workstation with two quad-core, hyper-threaded 2.93GHz Intel Xeon processors with 96GB of memory.  This yields a total of 16 threads allocated to each invocation of  CPLEX. The default behavior of CPLEX 12.1 is to allocate a number of threads equal to the number of machine cores. In the case of hyper-threaded architectures, each core is presented as a virtual dual-core -- although it is important to note that the performance is not equivalent to a true dual core. The workstation is shared by other users, such that our run-time results should be interpreted as conservative. With the exception of the optimality gap, which we set to 0.1\%, we used the CPLEX default settings for all other parameters. All runs were allocated a maximum of $10,800$ seconds (3 hours) of wall clock time.

We executed our models and algorithms on five test systems of varying size: the 6-bus, IEEE 24-bus, RTS-96, and IEEE 118-bus test systems \cite{IEEEtest}, and a simplified model of the US Western interconnection (WECC-240)\cite{Price2011}. The 6-bus system described in Section \ref{ieee6_example} is further augmented with three fast ramping generation units located at bus 1, 2, and 6, respectively, to ensure there is sufficient generation capacity for larger-size contingencies.  Generator data for these three units are identical to G4-G6 in Table \ref{tab10}.  To ensure there is sufficient operational flexibility in the WECC-240  system, we made eight transmission lines and one generation unit immune to failures. These nine elements include serial lines, pairs of transmission lines, and generation unit and transmission line pairs, whose failure would result in islanding of subsystems (buses). Additionally, we assume that non-dispatchable generation injections into the system can be shed during contingency states. For each test system, we consider a 24 hour planning horizon and the four contingency budgets $k=0, 1, 2,$ and $3$, yielding a total of 20 instances.

 \begin{table}[t]
\caption{Runtimes for different solution approaches to the CCUC problem}
\centering
\small{
\begin{tabular}{l c c c| c c c}
\hline
\multicolumn{4}{c|}{} & \multicolumn{3}{c}{Solution time} \\[-0.3ex]
\multicolumn{4}{c|}{} & \multicolumn{3}{c}{(exit status or feasibility gap)} \\[0.5ex]
Test System 		&$C$ 		& $k$	&  $\varepsilon_ {k}$ 		& EF        				 & BD        		& CSA2           	 	\\[0.5ex]
\hline

 6-bus			&	0			&	0	&	0 		    					& 0 				& 0 			 & 0			 	\\[-0.3ex]
         				&	16			&	1	&	0 		    					& 3 				& 2 			 & 1				\\[-0.3ex]
           			&	136			&	2	&	0.29 							& 7 				& 16 			 & 2				 \\[-0.3ex]
        				&	696			&	3	&	0.77 							& 134 			& 189 		 & 4				 \\[-0.3ex]
\hline
IEEE 24-bus     		&	0			&	0	&	0 		    					& 0 				& 1 			& 0				\\[-0.3ex]
         				&	70			&	1	&	0 		    					& 108 			& 58			& 11				\\[-0.3ex]
            			&	2,485		&	2	&	0.08 							& x(LPR) 			& 3,861 		& 101			\\[-0.3ex]
        				&	57,225		&	3	&	0.21 							& x(OM)			& x(0.03) 		& 397			 \\[-0.3ex]
\hline
RTS-96     		&	0			&	0	&	0 		    					& 1 				& 2 			 & 2				\\[-0.3ex]
         				&	216			&	1	&	0 		    					& x(LPR)			& 303 		& 41				\\[-0.3ex]
            			&	23,436		&	2	&	0.05 							& x(OM) 			& 8,139 	 	& 4,04			\\[-0.3ex]
        				&	1,679,796		&	3	&	0.09 							& x(OM) 			 & x(0.050) 		& 4,989			 \\[-0.3ex]
\hline
IEEE 118-bus     		&	0		&	0	&	0 		    					& 1 				& 1 			& 1				\\[-0.3ex]
         				&	240			&	1	&	0.01 		    					& x(LPR) 			 & 3,513 		& 352			\\[-0.3ex]
            			&	28,920		&	2	&	0.12 							& x(OM) 			& x(0.204) 		& 1,232			\\[-0.3ex]
        				&	2,304,200		&	3	&	0.25 							& x(OM) 			 & x(0.249) 		& 8,586			 \\[-0.3ex]
\hline
WECC-240     		&	0			&	0	&	0 		    					& 1 				& 2 			 & 2				\\[-0.3ex]
         				&	424			&	1	&	0 		    					& x(LPR) 			 & 262 		& 108			\\[-0.3ex]
            			&	90,100		&	2	&	0.06 							& x(OM) 			& x(0.004) 		& 2,484			\\[-0.3ex]
        				&	12,704,524	&	3	&	0.15 							& x(OM) 			& x(0.004) 		& x(0)			\\[-0ex]
\hline
\end{tabular}}
\label{tab1}
\end{table}

We first consider the run-times for the three different algorithms for solving the CCUC problem: the extensive form MILP, Benders decomposition, and the Contingency Screening Algorithm 2 (CSA2). The results are presented in Table \ref{tab1}. All times are reported in wall clock (elapsed) seconds. The column labeled ``$C$" reports the number of distinct contingencies for a given budget $k$, while the column labeled ``$\varepsilon_k$" reports the fraction of total load (demand) that can be shed. Entries in Table \ref{tab1} reporting ``x" indicate that the corresponding algorithm failed to locate a $N$-$k$-$\boldsymbol \varepsilon$ compliant solution within the 0.1\% optimality gap within the 3 hour time limit.  For those instances that could not be solved within the allocated time, we provide exit status or feasibility gaps, indicating the maximum fraction of total demand shed \emph{above} the allowable threshold $\varepsilon_k$ in the final solution. In all runs of the CSA2 algorithm, we initialize the contingency list $\mathscr L$ to the empty list.

As expected, the extensive form approach (EF) can only solve the smallest  instances, since for each contingency, a full set of DC power flow constraints (\ref{cont_const}) must be explicitly embedded in the formulation. As the number of contingencies grows, this formulation quickly becomes intractable.  The exit status ``LPR" and ``OM'' represent ``solving Linear Programming Relaxation at root node'' and ``Out of Memory'',  respectively. Note that our test instances only represent small to at most moderate sized systems relative to real power systems (which can contain on the order of thousands to tens of thousands of elements), indicating that even significant advances in solver technology are unlikely to mitigate this issue. Further, even given significant algorithmic advances, the memory requirements associated with the EF will likely cause the intractability to persist.

The BD approach attempts to address the exponential, as shown in Remark 6, explosion in the number of contingencies via a Benders reformulation/decomposition, with corresponding delayed cut generation. However, although the BD approach does not explicitly incorporate  power flow constraints (\ref{cont_const}) for each contingency into the formulation, those power flow constraints  must still be solved to identify violated feasibility cuts (which are then added to the master problem). In summary, the BD approach mitigates the memory issues associated with the EF approach, but the cost of identifying feasibility violations for a rapidly growing number of contingencies remains prohibitive. Overall, the BD approach can solve larger instances than the EF approach, but still fails given larger $k$ and larger test instances.

Finally, we consider the performance of our third approach: CSA2. Here, we see that all of our test instances, with the sole exception of the WECC-240 system with $k=3$, can be solved within the 3 hour time limit. This result is enabled by the combination of using a dynamic contingency list (significantly reducing the number of PSIP solves) and the fact that we are able to implicitly evaluate all the contingencies in order to identify a violated contingency, and then quickly find a corresponding feasibility cut by solving a single linear program (DCF). These features of the CSA2 algorithm allow it to mitigate the effects of a combinatorial number of contingencies and the associated impact on run-times and memory requirements.  Lastly, we note that although CPA2 failed to solve the WECC-240 system with $k=3$ within the allocate time, the final solution at the three hour mark is in fact a $N$-$k$-$\boldsymbol \varepsilon$ compliant solution.  For large power systems and/or contingency budgets, significant computational time is required to ``prove'' feasibility.  Eliminating the three hour time limit, we observed that the WECC-240 system with $k=3$ could be solved in approximately 18 hours, with the majority of this time taken to prove feasibility of the final solution.

We next examine the run-times of our CSA2 algorithm in further detail, as reported in Table \ref{tab2}. For each instance, we report the total number of possible contingency states $C$ and the number of contingencies for which corresponding feasibility cuts were actually generated. The latter corresponds to the final size of the dynamic contingency list, which is reported in the column labeled ``$|\mathscr L|$". Clearly, $|\mathscr L|$ corresponds to a vanishingly small fraction of the possible number of contingencies, which is critical to the tractability of the approach. The remaining columns of Table \ref{tab2} break down the total run time (in wall clock seconds) by the three main components of the algorithm -- the RMP, which identifies unit commitments; the power system inhibition problem (PSIP), which identifies a contingency that has no feasible corrective recourse power flow given the current RMP UC decisions and no-contingency economic dispatch; and the contingency feasibility subproblems (DCF), which yield the feasibility cuts.  The final column, labeled ``cuts'', reports the total number of feasibility cuts generated in solving the instance. It is clear from Table \ref{tab2} that the computational bottleneck in the CSA2 algorithm is the solution of the PSIP, such that any improvements to that process will yield immediate reductions in CSA2 run-times.

\begin{table}[t]
\caption{Runtime breakout for the CSA2 algorithm}
\centering
\footnotesize{
\begin{tabular}{l r r r | r r r | r r r r }
\hline
Test Systems 		&$C$ 		& 	$k$	&  $\varepsilon_ k$ 		&RMP        	&PSIP       		 &DCF        		& itr		&$|\mathscr L|$ 	&cuts 	\\[0.5ex]
\hline
 6-bus			&	0				&	0	&	0 		    		& 0 		    	& 0		        	 & 0				& 1		&0		&0			\\[-0.3ex]
         				&	16				&	1	&	0 		    		& 0 			& 1 			 & 0				& 2		&1		&21			\\[-0.3ex]
        				&	136				&	2	&	0.29 				& 0 			& 2 			 & 0			    	& 7		&3		&48			\\[-0.3ex]
        				&	696				&	3	&	0.77 				& 0 			& 4 			 & 0			    	& 11		&4		&89			\\ [-0.3ex]
\hline
IEEE 24-bus    		&	0				&	0	&	0 		    		& 0 		    	& 0		        	 & 0				& 1		&0		&0			\\[-0.3ex]
       				&	70				&	1	&	0 		    		& 6			& 46 			& 1				 & 185	&1		&185			\\[-0.3ex]
        				&	2,485			&	2	&	0.08				& 26 			& 69 			 & 6			   	& 857	&3		&857			\\[-0.3ex]
        				&	57,225			&	3	&	0.21 				& 64 			& 324 		& 9			   	 & 928	&4		&928			\\ [-0.3ex]
\hline
RTS-96    			&	0				&	0	&	0 		    		&2 		    	& 0		        	 & 0				& 1		&0		&0			\\[-0.3ex]
       				&	216				&	1	&	0 		    		& 13			& 25 			& 3				 & 12		&3		&27			\\[-0.3ex]
        				&	23,436			&	2	&	0.05				& 17			& 385 		& 2			   	 & 15		&4		&33			\\[-0.3ex]
        				&	1,679,796			&	3	&	0.09 				& 19 			& 4,965 		 & 5			   	& 17		&5		&38			\\ [-0.3ex]
\hline
IEEE 118-bus    	&	0				&	0	&	0 		    		&1 		    	& 0		        	 & 0				& 1		&0		&0			\\[-0.3ex]
       				&	240				&	1	&	0.01 		    		& 243		& 72 			& 37				 & 85		&5		&1,305			\\[-0.3ex]
        				&	28,920			&	2	&	0.12				& 377 		& 796 		& 59				 & 120	&7		&1,671		\\[-0.3ex]
        				&	2,304,200			&	3	&	0.25 				& 405 		& 8,114 		 & 67				& 132	&9		&1,743		\\ [-0.3ex]
\hline
WECC 240-bus    		&	0			&	0	&	0 		    		&0 		    	& 0		        	 & 0				& 1		&0		&0			\\[-0.3ex]
       				&	424				&	1	&	0 		    		& 4			& 102 		& 2				 & 5		&2		&48			\\[-0.3ex]
        				&	90,100			&	2	&	0.06				& 3 			&2,479 		 & 2				& 5		&2		&48		\\[-0.3ex]
        				&	12,704,524		&	3	&	0.15 				& x 			& x 			 & x				& x		&x		&x			\\ [-0ex]
\hline
\end{tabular}}
\label{tab2}
\end{table}

\section{Conclusions}
\label{sec:conclusions}

We have investigated the problem of committing generation units in power system operations, and determining a corresponding no-contingency state economic dispatch, such that the resulting solution satisfies the $N$-$k$-$\boldsymbol \varepsilon$ reliability criterion. This reliability criterion is a generalization of the well-known $N$-$k$  criterion, and requires that at least $(1-\varepsilon_ j)$ fraction of the total demand is met following the failure of $j$  system components, for $j \in \{1,\cdots,k\}$. We refer to this problem as the contingency-constrained unit commitment problem, or CCUC. We proposed two algorithms to solve the CCUC: one based on the Benders decomposition approach, and another based on contingency screening algorithms. The latter method avoids the combinatorial explosion in the number of contingencies by seeking vulnerabilities in the current solution, and generating valid inequalities  to exclude such infeasible solutions in the master problem.  We tested our proposed algorithms on test systems of varying sizes.  Computational results show our proposed Contingency Screening Algorithm (CSA2), which uses a bi-level separation program to implicitly consider all contingencies and a dynamic contingency list to avoid re-identification of contingencies, significantly outperforms the Benders decomposition approach.  We were able to solve all test systems, with the exception of the largest WECC-240 instance, in under 3 hours. In contrast, both the Benders decomposition algorithm and a straightforward solution of the CCUC extensive form, failed to solve all but the smallest instances within 3 hours.

We believe that  this paper will provide a significant basis for subsequent research in contingency-constrained unit commitment. For example, we are working to apply these methods to full-scale systems. While our results are promising in terms of scalability, full-scale problems pose more significant computational challenges, and consequently will require stronger formulations for the power system inhibition problem and possible adoption of high-performance computing resources.  Further, our current CCUC model assumes  all  component failures occur simultaneously. In order to reflect practical operational situations, where failures may happen consecutively, new CCUC models that consider timing between system component failures are needed. We plan to extend our CCUC models to include these cases. Finally, we worked exclusively with a deterministic CCUC model to date. However, it is ultimately essential to take uncertainty into account in unit commitment, e.g., to account for uncertain demand and renewable generation units. We believe our current cutting plane framework can be naturally extended to robust optimization and stochastic programming formulations via a nested decomposition approach.

\vskip 0.5cm

{\bf\it Acknowledgement}. {
Sandia National Laboratories' Laboratory-Directed Research and Development Program and the U.S. Department of Energy's Office of Science (Advanced Scientific Computing Research program) funded portions of this work. Sandia National Laboratories is a multi-program laboratory managed and operated by Sandia Corporation, a wholly owned subsidiary of Lockheed Martin Corporation, for the U.S. Department of Energy's National Nuclear Security Administration under contract DE-AC04-94AL85000.}

\bibliographystyle{}

\begin{thebibliography}{0}

\bibitem{Arroyo2010}%[Arroyo (2010)]
Arroyo, J.M. 2010. Bilevel programming applied to power system vulnerability analysis under multiple contingencies. {\em IET Gener. Transm. Distrib.} {\bf 4}(2): 178--190.

\bibitem{Benders1962}%[Benders (1962)]
Benders, J. 1962. Partitioning procedures for solving mixed-variables programming problems. {\em Numerische Mathematik} {\bf 4}: 238--252.

\bibitem{Bienstock2010}%[Bienstock and Verma (2010)]
Bienstock, D., A. Verma. 2010. The N-k problem in power grids: new models, formulations, and numerical experiments. {\em SIAM J. Optim.} {\bf 20}(5): 2352--2380.

%\bibitem{Burch2003}%[Burch et al. (2003)]
%Burch, C., R. Carr, S. Krumke, M. Marathe, C. Phillips and E. undberg. 2003. A decomposition-based approximation for network inhibition. {\em Network Interdiction and Stochastic Integer Programming}, D.L. Woodruff, eds., pp. 51--66.

\bibitem{Capitanescu2011}%[Capitanescu et al. (2011)]
Capitanescu, F., et al. 2011. State-of-the-art, challenges, and future trends in security constrained optimal power flow. {\em Elec. Power Syst. Res.} {\bf 81}: 1731--1741.

%\bibitem{Caroe1998}%[Car{\o}e and Schultz (1998)]
%Car{\o}e, C.C., R. Schultz. 1998. A two-stage stochastic program for unit commitment under uncertainty in a hydro-thermal power system. {\em Preprint SC 98-11}, Konrad-Zuse-Zentrumf\"{u}r Informationstechnik Berlin.

\bibitem{Carrion2006}%[Carri\'{o}n and Arroyo (2006)]
Carri\'{o}n, M., J.M. Arroyo. 2006. A compuationally efficient mixed-integer linear formulation for the thermal unit commitment problem. {\em IEEE Trans. Power Syst.} {\bf 21}(3): 1371--1378.

%\bibitem{Cerisola2009}%[Cerisola et al. (2009)]
%Cerisola, S., A. Baillo, J.M. Fernandez-Lopez, A. Ramos, R. Gollmer. 2009. Stochastic power generation unit commitment in electricity markets: A novel formulation and a comparison of solution methods. {\em Oper. Res.} {\bf 57}(1) 32--46.

\bibitem{Chen2012a}%[Chen et al. (2012)]
Chen, R.L., A. Cohn, N. Fan, A. Pinar. 2012. $N$-$k$-$\epsilon$ power system design. {\em Proc. 12th Probabilistic Methods for Power Syst. Conf.}. Istanbul, Turkey.

\bibitem{Chen2014}
Chen, R.L., A. Cohn, N. Fan, A. Pinar. 2014. Contingency-risk informed power system design. {\em IEEE Trans. Power Systems}. DOI: 10.1109/TPWRS.2014.2301691.

%\bibitem{elakirby08}%[Ela and Kirby (2008)]
%Ela. E., B. Kirby. ERCOT Event on February 26, 2008: Lessons Learned. Technical Report, NREL/TP-500-43373, July 2008.

\bibitem{Fan2011}%[Fan et al. (2011)]
Fan, N., H. Xu, F. Pan, P.M. Pardalos. 2011. Economic analysis of the N-k power grid contingency selection and evaluation by graph algorithms and interdiction methods. {\em Energy Syst.} {\bf 2}(3-4): 313--324.

\bibitem{fercreport}%[FERC (2011)]
FERC Staff Report. Recent ISO Software Enhancements and Future Software and Modeling Plans. Federal Energy Regulatory Commission.

\bibitem{Fu2005}%[Fu et al. (2005)]
Fu, Y., M. Shahidehpour, Z. Li. 2005. Security-constrained unit commitment with AC constraints. {\em IEEE Trans. Power Syst.} {\bf 20}(3): 1538-1550.

\bibitem{Fu2006}%[Fu et al. (2006)]
Fu, Y., M. Shahidehpour, Z. Li. 2006. AC contingency dispatch based on security-constrained unit commitment. {\em IEEE Trans. Power Syst.} {\bf 21}(2): 897--908.

\bibitem{Hedman2010}%[Hedman et al. (2010)]
Hedman, K.W., M.C. Ferris, R.P. O'Neill, E.B. Fisher, S.S. Oren. 2010. Co-optimization of generation unit commitment and transmission switching with N-1 reliability. {\em IEEE Trans. Power Syst.} {\bf 24}(2): 1052--1063.

\bibitem{Hobbs2001}%[Hobbs et al. (2001)]
Hobbs, B. F., M. H. Rothkopf, R. P. O'Neill, H.-P. Chao, eds. 2001. {\em The Next Generation of Electric Power Unit Commitment Models}. Kluwer Academic Publishers, Boston.

\bibitem{IEEEtest}%[IEEE Test Systems]
IEEE Test Systems. \url{http://www.ee.washington.edu/research/pstca}.

\bibitem{LeRoDoPi06}   %[Lesieutre et al. ( 2006)]
B. Lesieutre, S. Roy, V. Donde, and A. Pinar, {Power system extreme event analysis
using graph partitioning}, {\em Proc.  39th North American Power Symposium}, Carbondale,
IL, October 2006.

\bibitem{LePiRo08} %[Lesieutre et al. (2008)]
B. Lesieutre, A. P{\i}nar, and S. Roy, Power System Extreme Event Detection: The
Vulnerability Frontier, in {\it Proc. 41st Hawaii International Conference on System Sciences}, pages 184, Waikoloa, Big Island, HI, 2008.

\bibitem{Lotfjou2010}%[Lotfjou et al. (2010)]
Lotfjou, A., M. Shahidehpour, Y. Fu, Z. Li. 2010. Security-constrained unit commitment with AC/DC transmission systems. {\em IEEE Trans. Power Syst.} {\bf 25}(1): 531--542.

\bibitem{NERC2011}%[NERC (2011)]
North American Electric Reliability Corporation, Transmission Planning Standards, Accessed on April 2014.
Available at http://www.nerc.com/pa/Stand/Reliability\%20Standards/Forms/AllItems.aspx
%North American Electric Reliability Corporation (NERC). Standard TPL-001-1--system performance under normal conditions. February 2011. Available at http://www.nerc.com/docs/standards/sar/Project\_2006-02\_TPL-001-1.pdf

\bibitem{ONeill2010}%[O'Neill et al. (2010)]
O'Neill, R.P., K.W. Hedman, E.R. Krall, A. Papavasiliou, S.S. Oren. 2010. Economic analysis of the N-1 reliable unit commitment and transmission swtiching problem using duality concepts. {\em Energy Syst.} {\bf 1})(2): 165--195.

\bibitem{personalcommunication}%[PersonalCommunication (2012)]
Personal Communication. Dr. Eugene Litvinov. March, 2012.

\bibitem{Padhy2004}%[Padhy (2004)]
Padhy, N.P. 2004. Unit commitment -- A bibliographical survey. {\em IEEE Trans. Power Syst.} {\bf 19}(3): 1196--1205.

\bibitem{Pinar2010}%[Pinar et al. (2010)]
Pinar, A., J. Meza, V. Donde, B. Lesieutre. 2010. Optimization strategies for the vulnerability analysis of the electric power grid. {\em SIAM J. Optim.} {\bf 20}(4): 1786--1810.

\bibitem{Price2011}%[Price (2011)]
Price, J.E. 2011. Reduced Network Modeling of WECC as a Market Design Prototype. {\em Proceedings of the 2011 IEEE General Meeting of the Power and Energy Society.} San Diego, California.

\bibitem{Salmeron2004}%[Salmeron et al. (2004)]
Salmeron, J., K. Wood, R. Baldick. 2004. Analysis of electric grid security under terrorist threat. {\em IEEE Trans. Power Syst.} {\bf 19}(2): 905--912.

\bibitem{Salmeron2009}%[Salmeron et al. (2009)]
Salmeron, J., K. Wood, R. Baldick. 2009. Worst-case interdiction analysis of large-scale electric power grids. {\em IEEE Trans. Power Syst.} {\bf 24}(1): 96--104.

%\bibitem{Shiina2004}%[Shiina and Birge (2004)]
%Shiina, T., J.R. Birge. 2004. Stochastic unit commitment problem. {\em Intl. Trans. in Op. Res.} {\bf 11} 19--32.

\bibitem{Street2011a}%[Street et al. (2011a)]
Street, A., F. Oliveira, J.M. Arroyo. 2011. Contingency-constrained unit commitment with n-K security criterion: A robust optimization approach. {\em IEEE Trans. Power Syst.} {\bf 26}(3): 1581--1590.

\bibitem{Street2011b}%[Street et al. (2011b)]
Street, A., F. Oliveira, J.M. Arroyo. 2011. Energy and reserve scheduling under an N-K security criterion via robust optimization. {\em Proc. 17th Power Syst. Compt. Conf.}. Stockholm, Sweden.

\bibitem{Wang2008}%[Wang et al. (2008)]
Wang, J., M. Shahidehpour, Z. Li. 2008. Security-constrained unit commitment with volatile wind power generation. {\em IEEE Trans. Power Syst.} {\bf 23}(3): 1319--1327.

\bibitem{Wang2012}%[Wang et al. (2012)]
Wang, Q., J.-P. Watson, Y. Guan. 2013. Two-stage robust optimization for $N$-$k$ contingency-constrained unit commitment. {\em IEEE Trans. Power Syst.} {\bf 28}(3): 2366--2375.

\bibitem{Wood1996}%[Wood and Wollenberg (1996)]
Wood, A.J., Wollenberg, B.J., 1996. {\em Power Generation, Operation and Control (2nd edition)}. John Wiley \& Sons, New York.

\bibitem{Wu2010}%[Wu and Shahidehpour (2010)]
Wu, L., M. Shahidehpour. 2010. Accelerating the Benders decomposition for network-constrained unit commitment problems. {\em Energy Syst.} {\bf 1}(3): 339--376.

\bibitem{Zeng2014}
Yuan, W., L. Zhao, B. Zeng. 2014. Optimal power grid protection through a defender-attacker-defender model. {\em Reliability Engineering and System Safety} {\bf 121}: 83--89.

\bibitem{Zeng2013}
Zhao, L., B. Zeng. 2013. Vulnerability analysis of power grids with line switching. {\em IEEE Transactions on Power Syst.} {\bf 28}(3): 2727--2736.

\bibitem{Zheng2013}
Zheng, Q.P., J. Wang, P.M. Pardalos, Y. Guan. 2013. A decomposition approach to the two-stage stochastic unit commitment problem. {\em Annals of Operations Research} {\bf 210}(4): 387--410.

\end{thebibliography}

\end{document}